\documentclass[journal]{IEEEtran}
\usepackage{etex}
\reserveinserts{100}
\usepackage{mathtools}
\usepackage{graphics, theorem, times, amsfonts, graphicx, amssymb, cite}
\usepackage{tikz}
\usepackage[T1]{fontenc}
\usetikzlibrary{shapes,arrows}
\usepackage{pgfplots}
\usepackage{color}
\usepackage{setspace}
\usepackage{hyperref}
\usepackage{multirow}
\usepackage{rotating}
\usepackage{comment}
\usepackage{flushend}
\usepackage[font=small]{caption}
\usepackage{subcaption}
\usepackage[affil-it]{authblk}
\usepackage{bm}
\usepackage{algorithm}
\usepackage{algorithmic}
\usepackage{enumerate}
\usepackage{morefloats}
\input{mysymbol.sty}
\usepackage{needspace}

\newcommand{\QED}{\hfill\ensuremath{\blacksquare}}

\newtheorem{assumption}{\hspace{0pt}\bf Assumption}
\newtheorem{lemma}{\hspace{0pt}\bf Lemma}

\newtheorem{theorem}{\hspace{0pt}\bf Theorem}
\newtheorem{corollary}{\hspace{0pt}\bf Corollary}

\newtheorem{remark}{\hspace{0pt}\bf Remark}

\def\df{{\nabla f}}

\date{\today}

\title{A Primal-Dual Quasi-Newton Method for Exact Consensus Optimization}
\author{Mark Eisen, Aryan Mokhtari, and Alejandro Ribeiro 
\thanks{{Supported by ARL DCIST CRA W911NF-17-2-0181 and Intel Science and Technology Center for Wireless Autonomous Systems (ISTC-WAS). M. Eisen and A. Ribeiro are with the Department of Electrical and Systems Engineering, University of Pennsylvania. A. Mokhtari is with the Laboratory for Information and Decision Systems, Massachusetts Institute of Technology. Email at: maeisen@seas.upenn.edu, aryanm@mit.edu, aribeiro@seas.upenn.edu. Part of the results in this paper appeared in \cite{eisen2017primal} . {This paper expands the results and presents convergence proofs that are referenced in \cite{eisen2017primal}.} 
}}}

\begin{document}
\thispagestyle{empty}
\maketitle

\begin{abstract}
We introduce the primal-dual quasi-Newton (PD-QN) method as an approximated second order method for solving decentralized optimization problems. The PD-QN method performs quasi-Newton updates on both the primal and dual variables of the consensus optimization problem to find the optimal point of the augmented Lagrangian. By optimizing the augmented Lagrangian, the PD-QN method is able to find the exact solution to the consensus problem with a linear rate of convergence. We derive fully decentralized quasi-Newton updates that approximate second order information to reduce the computational burden relative to dual methods and to make the method more robust in ill-conditioned problems relative to first order methods. The linear convergence rate of PD-QN is established formally and strong performance advantages relative to existing dual and primal-dual methods is shown numerically.
\end{abstract}


\section{Introduction} \label{sec_intro}

In the setting of decentralized optimization, a connected network of agents, or nodes, are interested in minimizing a common global objective function, the components of which are distributed locally across all the agents. To jointly optimize the global objective, nodes must collaborate with their neighbors by successively sharing information of their locally measurable objective function components. Decentralized optimization has proven effective in contexts where information is gathered by different nodes of a network, such as decentralized control
\cite{Bullo2009,Cao2013-TII,LopesEtal8}, wireless systems
\cite{Ribeiro10,scutari2014decomposition}, sensor networks
\cite{Schizas2008-1,KhanEtal10,cRabbatNowak04}, and large scale
machine learning
\cite{bekkerman2011scaling,Cevher2014}.

Perhaps the most common and well-studied problem in decentralized optimization is the consensus optimization problem. Here, the common minimizer that optimizes the global objective is found locally at each node through the use of local copies of the decision variable. Each node minimizes its local objective using its local copy while simultaneously seeking agreement, or \emph{consensus}, with its neighbors. This approach allows for a fundamentally decentralized manner of optimizing a global function that relies only on the ability to exchange information with neighbors. There are many methods that solve the consensus problem, differing largely in how the consensus condition is enforced. Among the popular techniques include the use of additional penalties for violating consensus in the objective function \cite{nedic2009, YuanQing, Jakovetic2014-1,shi2015extra,mokhtari2016dsa,mokhtari2017network}. Alternatively, the consensus can formulated as an explicit constraint, which can be optimized directly in the dual domain {\cite{cRabbatNowak04,Schizas2008-1,makhdoumi2016convergence,bianchi2014stochastic,zargham2014accelerated}}. A promising extension of dual-based methods include ``primal-dual'' methods, which iteratively find solutions that incrementally get closer to both optimality and consensus \cite{chang2014distributed,mokhtari2016decentralized,li2017primal}. These methods are beneficial in combining the low computational cost of penalty-based methods with the exactness of dual-based methods.

In general, standard consensus optimization techniques that rely only on first order information contained in the gradient, such as gradient descent, suffer from slow convergence rates. These slow convergence rates are particularly apparent in problems that are ill-conditioned, or in other words have large spread in eigenvalues of the function's Hessian matrix. In centralized optimization, the ill-conditioning is commonly corrected by incorporating second order Hessian information into the descent computation using the Newton step.  While this technique cannot be used directly in distributed optimization due to the non-sparsity of the Hessian inverse, there exist ways of using second order information to approximate the Newton step for distributed settings. This has been done for consensus optimization problems reformulated as both the penalty-based methods \cite{mokhtari2017network,mansoori2017superlinearly}  and dual-based methods \cite{zargham2014accelerated}, as well as the more recent primal-dual methods \cite{mokhtari2015dqm,mokhtari2016decentralized}. These approximate Newton methods exhibit faster convergence relative to their corresponding first order methods.

Despite the advances in distributed Newton-based methods, there are many cases in which the exact Hessian information is either difficult or computationally expensive to compute. The centralized alternative to Newton's method comes in the form of quasi-Newton methods, which use gradients to produce a curvature estimation in lieu of the Hessian inverse \cite{dennis1974characterization,powell1976some}. In the distributed setting, the commonly used Broyden-Fletcher-Goldfarb-Shanno (BFGS) quasi-Newton has been adapted for both penalty-based and dual-based formulations \cite{eisen2017decentralized}. This method improves convergence but suffers from many of the same issues of first-order penalty-based and dual-based methods regarding accuracy and computational burden. In this paper, we develop a novel ``primal-dual'' quasi-Newton method that achieves state of the art convergence rates to exact solutions at a lower computational cost relative to the dual-based alternative. This is done through the approximation of an internal optimization problem present in dual methods with a fully distributed quasi-Newton update. We further employ a distributed dual quasi-Newton update to accelerate the dual ascent relative to dual gradient updates. In this way, we can obtain a method with a linear convergence rate to the exact solution with the computational complexity of penalty-based methods, all while not requiring the computation of Hessian information.  

The paper begins with a formal statement of the consensus optimization problem and the introduction of the augmented Lagrangian (Section \ref{sec_problem_formulation}). For solving the consensus problem in a distributed manner, a set of standard dual methods (Section \ref{sec_dual_methods}) and augmented Lagrangian methods (Section \ref{sec_aug_lag}) are discussed. The former methods suffer from slow convergence rates and the latter methods suffer from either the inability to distribute the computations or the necessity in solving an internal minimization problem at every step. All of these issues are bypassed by substituting the standard primal and dual updates in these dual methods with quasi-Newton updates. This provides the basis of the proposed PD-QN method (Section \ref{sec_dbfgs}). Separate quasi-Newton updates are derived for both the primal update (Section \ref{sec_primal_update}) and dual update (Section \ref{sec_dual_update}), each of which is designed to be distributedly computable while retaining desirable properties of traditional quasi-Newton methods. Convergence properties are then established (Section \ref{sec_convergence}). Given standard properties of smoothness and strong convexity, we demonstrate a linear convergence rate to the exact solution of the consensus problem. We close the paper with numerical results comparing the performance of PD-QN to first and second order methods on various consensus problems of practical interest (Section \ref{sec_numerical_results}).

\section{Problem Formulation} \label{sec_problem_formulation}

We consider a distributed system of $n$ nodes connected by an undirected communication graph $\ccalG=(\ccalV,\ccalE)$ with nodes $\ccalV=\{1,\dots,n\}$ and $m$ edges $\ccalE=\{(i,j)\ |\ i\ \text{and}\ j \ \text{are connected} \}$. Define the set $n_i$ as the neighborhood of node $i$ including $i$, i.e., $n_i=\{j\ |\ j=i \lor (i,j)\in\ccalE\}$, and the neighborhood size $m_i := |n_i|$. Each node $i$ has access to a local strongly convex cost function $f_i: \reals^p \rightarrow \reals$ and the goal is to find the optimal variable $\tbx^* \in \reals^p$ that minimizes the aggregate of all local cost functions $\sum_i f_i$ with distributed computations. By distributed computations, we mean in particular that each node is able to itself obtain the common minimizer $\tbx^*$ through computations performed locally and through information exchanges with its neighbors in $\ccalE$. To solve this locally, we consider the consensus formulation, in which each node $i$ stores and maintains a local variable $\bbx_i \in \reals^p$ and seeks to minimize its local cost $f_i(\bbx_i)$ while satisfying a consensus constraint with all neighbors. More specifically, consider the global variable $\bbx=[\bbx_1;\dots;\bbx_n]\in\reals^{np}$ and resulting optimization problem of interest
\begin{align}\label{eq_primal_problem}
   \bbx^* \ &:= \ \argmin_{\bbx\in \reals^{np}} \ f(\bbx) = \sum_{i=1}^n f_i(\bbx_i)  \\
               &\quad \st \  (\bbI-\bbZ)^{1/2}\bbx=\bb0, \nonumber
\end{align}
where the matrix $\bbZ\in \reals^{np \times np}$ is chosen so that the feasible variables in \eqref{eq_primal_problem} satisfy the consensus constraint $\bbx_i=\bbx_j$ for all $i,j$. A customary choice of a matrix $\bbZ$ with this property is to make it the Kronecker product $\bbZ \coloneqq \bbW \otimes \bbI_p$ of a weight matrix $\bbW\in\reals^{n\times n}$ and the identity matrix $\bbI_p\in\reals^{p \times p}$. The elements of the weight matrix are $w_{ij}> 0$ if $(i,j)\in\ccalE$ and $w_{ij}= 0$ otherwise and the weight matrix $\bbW\in \reals^{n\times n}$ is further assumed to satisfy
\begin{equation}\label{weight_matrix_conditions}
   \bbW                     = \bbW^T,      \quad 
   \bbW\mathbf{1}           = \mathbf{1},  \quad 
   \text{null}\{\bbI-\bbW\} = \text{span}\{\mathbf{1}\}.
\end{equation}
The first two conditions in \eqref{weight_matrix_conditions} imply symmetry and row stochasticity of $\bbW$, respectively. The third condition enforces the consensus constraint. We further impose the following assumption on diagonal weights necessary for the analysis later in this paper.
 \begin{assumption} \label{as_weight_bound}
The diagonal entries of $\bbW$ are bounded by two positive constants $1>\Delta > \delta>0$. I.e., for all $i$ we have that $\delta<w_{ii}<\Delta$.
\end{assumption} 
Since $\text{null}(\bbI-\bbW)=\text{null}(\bbI-\bbW)^{1/2}=\text{span}\{\mathbf{1}\}$, it follows that for any vector $\bbx=[\bbx_1;\dots;\bbx_n]\in \reals^{np}$ the relation $(\bbI-\bbZ)^{1/2}\bbx=\bb0$ holds if and only if $\bbx_1=\dots=\bbx_n$. This means that the feasible variables in \eqref{eq_primal_problem} indeed satisfy $\bbx_i=\bbx_j$ for all $i,j$ and that, consequently, the problem in \eqref{eq_primal_problem} is equivalent to the minimization of $\sum_{i=1}^n f_i(\tbx)$ and subsequently $\bbx^* = [\tbx^*; \tbx^*; \hdots; \tbx^*]$. 

To solve \eqref{eq_primal_problem}, it is necessary to form the Lagrangian. Define $\bbnu:= [\bbnu_1, \bbnu_2, \hdots, \bbnu_n] \in \reals^{np}$ to be a set of dual variables for the equality constraint $(\bbI-\bbZ)^{1/2}\bbx=\bb0$, with node $i$ holding the $i$th block $\bbnu_i \in \reals^p$. The Lagrangian is then defined as $\ccalL_0(\bbx,\bbnu) := \sum_{i=1}^n f_i(\bbx_i) + \bbnu^T (\bbI-\bbZ)^{1/2}\bbx$. To remove the dependence of the potentially dense matrix $(\bbI-\bbZ)^{1/2}$, we introduce an adjusted dual variable $\bby := (\bbI-\bbZ)^{1/2}\bbnu$ and work directly with $\bby \in \reals^{np}$. We additionally may add a quadratic penalty term to form the augmented Lagrangian in terms of $\bbx$ and $\bby$ as
\begin{align}\label{eq_aug_lagrangian}
\ccalL_{\alpha}(\bbx,\bby) := \sum_{i=1}^n f_i(\bbx_i) + \bby^T\bbx + \frac{\alpha}{2}\bbx^T(\bbI-\bbZ)\bbx,
\end{align}
indexed by $\alpha\geq0$, the weight of the quadratic penalty. By setting $\alpha=0$, we remove the quadratic penalty and recover the standard Lagrangian $\ccalL_0(\bbx,\bby)$. Note that any feasible point to the constrained problem in \eqref{eq_primal_problem} will make the quadratic penalty in \eqref{eq_aug_lagrangian} zero and thus add no additional cost to the standard Lagrangian. A pair $\bbx^*$ and $\bby^*$ that jointly optimize $\ccalL_{\alpha}(\bbx,\bby)$ will thus provide the solution to the original problem of interest in \eqref{eq_primal_problem} for any choice of $\alpha$.

\subsection{Dual methods} \label{sec_dual_methods}

A standard set of approaches towards jointly optimizing \eqref{eq_aug_lagrangian} is to operate exclusively on the dual variable $\bby$ in what are known as dual methods. The basic dual objective function to be maximized is obtained by minimizing a \emph{non-augmented} Lagrangian $\ccalL_{0}(\bbx,\bby)$ over $\bbx$, where $\ccalL_0(\bbx,\bby)$ is defined as the Lagrangian in \eqref{eq_aug_lagrangian} with no quadratic weight, i.e. $\alpha=0$. This results in the following dual problem
\begin{align}\label{eq_dual_problem}
\bby^* := \argmax_{\bby \in \reals^{np}} \left[\min_{\bbx \in \reals^{np}} \sum_{i=1}^n f_i(\bbx_i) + \bby^T \bbx \right].
\end{align}
The optimal primal variable given any dual iterate $\bby$ found in the internal minimization step in \eqref{eq_dual_problem} we notate as $\bbx_0(\bby)$. The optimal $\bbx^*$ from the original problem in \eqref{eq_primal_problem} can then be recovered as $\bbx^* = \bbx_0(\bby^*)$ due to strong duality of a strongly convex problem---see, e.g., \cite{boyd2004convex}. A wide array of optimization techniques can be used to iteratively solve for $\bby^*$. Gradient ascent, for example, can be performed directly on the maximization problem in \eqref{eq_dual_problem}. In gradient ascent, at each iteration index $t=0,1,2,\hdots$ the dual variable $\bby_{t+1}$ is computed using the previous iterate $\bby_{t}$ and the gradient of objective function in \eqref{eq_dual_problem}, resulting in the respective primal and dual updates,
\begin{align}
\bbx_{t+1} &= \argmin_{\bbx} \ccalL_{0}(\bbx,\bby_{t}) \label{eq_dual_ascent_p} \\
\bby_{t+1} &= \bby_{t} + \eps_d (\bbI - \bbZ) \bbx_{t+1}, \label{eq_dual_ascent_d}
\end{align}
where $\eps_d > 0$ is a scalar step size. Observe that, in both the updates in \eqref{eq_dual_ascent_p}-\eqref{eq_dual_ascent_d}, the $i$th block of $\bbx_{i,t+1}$ and $\bby_{i,t+1}$ can be computed locally by node $i$ using only local exchanges with neighbors $j \in n_i$. Further observe this distributed capability is permitted only by considering the non-augmented Lagrangian in the definition of the dual function in \eqref{eq_dual_problem}. The updates together are commonly known as dual ascent (DA) \cite{cRabbatNowak04} and are known to converge sub-linearly to the optimal pair $(\bbx^*, \bby^*)$ with a rate of $\ccalO(1/t)$. However, first order gradient-based methods can be further slowed in practice when the problem is ill-conditioned, motivating the development of more sophisticated dual techniques.

Quasi-Newton methods are a well known alternative to first order methods. In traditional centralized settings, they are known to have more desirable convergence properties and perform better in practice than first order methods because they approximate a curvature correction to the descent direction. They have recently been adapted for distributed algorithms as well.  In \cite{eisen2016decentralized,eisen2017decentralized}, a dual quasi-Newton method is derived that replaces the first order dual update in \eqref{eq_dual_ascent_d} with an approximate second order update of the form 
\begin{align}
\bby_{t+1} &= \bby_{t} +\eps_d \bbH^{-1}_{t}(\bbI - \bbZ) \bbx_{t+1}, \label{eq_dual_dbfgs}
\end{align}
where $\bbH_{t}$ is an approximation of the dual Hessian and thus serves as a distributed approximation to Newton's method. In particular, $\bbH_{t}$ is constructed as a positive definite matrix that satisfies what is known as the secant condition. Recall $\bbh_{t} = (\bbI - \bbZ)\bbx_{t+1}$ as the dual gradient, and define the dual variable variation $\bbv_{t}$ and gradient variation $\bbs_{t}$ vectors, 
\begin{align}
\bbv_{t} = \bby_{t+1} - \bby_{t}, \enskip \bbs_{t} = (\bbI - \bbZ) (\bbx_{t+1}-\bbx_{t}).\label{eq_d_bfgs_vars} 
\end{align}
Observe that $\bbv_{t}$ and $\bbs_{t}$ capture differences of two consecutive dual variables and gradients, respectively, evaluated at steps $t+1$ and $t$. At each iteration, we select a new Hessian approximation $\bbH_{t+1}$ that satisfies the secant condition $\bbH_{t+1} \bbv_{t} =  \bbs_{t}$. This condition is fundamental, as the secant condition is satisfied by the actual Hessian for small $\bbv_{t}$. The dual D-BFGS \cite{eisen2016decentralized} method designs a matrix $\bbH_{t}$ that satisfies the global secant conditions while being locally computable by each node. The inclusion of this approximation is shown to numerically improve upon the first order DA method, but ultimately suffers from the same slow converge rate of $\ccalO(1/t)$. 

\subsection{Augmented Lagrangian methods}\label{sec_aug_lag}
A well studied technique to improve upon the convergence rate of dual methods is to operate on the augmented Lagrangian $\ccalL_{\alpha}(\bbx,\bby)$ in \eqref{eq_aug_lagrangian} for some $\alpha >0$. Consider substituting the primal update in \eqref{eq_dual_ascent_p} with the minimization over the augmented Lagrangian
\begin{align}
\bbx_{t+1} &= \argmin_{\bbx} \ccalL_{\alpha}(\bbx,\bby_{t}). \label{eq_mm_p}
\end{align}
Using the augmented primal update in \eqref{eq_mm_p} with the first order dual update in \eqref{eq_dual_ascent_d} results in a method commonly referred to as the method of multipliers (MM), which is shown to exhibit a fast linear convergence rate to the optimal primal-dual pair \cite{bertsekas2014constrained,jakovetic2015linear,mokhtari2016decentralized}. However, MM cannot be implemented in a distributed manner because \eqref{eq_mm_p} cannot be computed locally due to the quadratic coupling term in \eqref{eq_aug_lagrangian}. The ADMM methods exists as a distributed alternative to MM, which permits distributed computation through a decoupling of the variables \cite{Shi2014-ADMM}. Both the augmented Lagrangian update in \eqref{eq_mm_p} and the standard update in \eqref{eq_dual_ascent_p} will in any case still require an internal minimization step. For most objective functions, this primal update will be computationally expensive, thus making these methods difficult to implement in practice.

Despite this shortcoming, it would nonetheless be beneficial to incorporate the augmented Lagrangian in the primal update to improve upon the convergence rate of distributed quasi-Newton methods. In this paper, we develop a primal-dual quasi-Newton method that, in addition to using the quasi-Newton dual update in \eqref{eq_dual_dbfgs} for improved conditioning, implements a quasi-Newton approximation of \eqref{eq_mm_p} that permits local and distributed computation to find exact solutions to \eqref{eq_primal_problem}.

\subsection{Related Work}

\begin{table*}[t]
\centering
\caption{Comparison of consensus optimization methods}
\label{tab_methods}
\begin{tabular}{l | lllll}
                        		& Convergence rate &  Exact solution? & Internal opt. problem? & Primal Order & Dual Order    \\ \hline
DGD \cite{nedic2009}          		& Linear   & N & N & First & N/A\\
Network Newton\cite{mokhtari2017network}        & Linear-Quadratic   & N & N & Second & N/A\\
D-BFGS (primal) \cite{eisen2017decentralized}        & Linear   & N & N & ``Second'' & N/A\\\hline
Dual descent  \cite{cRabbatNowak04}          & Sublinear   & Y & Y & N/A & First\\
ADMM \cite{Schizas2008-1,Shi2014-ADMM}            	& Linear   & Y & Y & N/A & First\\
D-BFGS (dual) \cite{eisen2017decentralized}         & Sublinear   & Y & Y & N/A & ``Second''\\\hline
EXTRA  \cite{shi2015extra,li2017primal}              	& Linear   & Y & N & First & First\\
ESOM  \cite{mokhtari2016decentralized}              	& Linear   & Y & N & Second & First\\
\textbf{PD-QN} 	& \textbf{Linear}   & \textbf{Y} & \textbf{N} & \textbf{``Second''} & \textbf{``Second''} \\
\end{tabular}
\end{table*}

The existing literature in solving the consensus problem in \eqref{eq_primal_problem} differs in many ways, ranging from convergence rate to computational complexity to communication cost. We summarize many of these methods in Table \ref{tab_methods} and discuss their qualities here. Both Table \ref{tab_methods} and the following discussion break consensus optimization methods into three basic classes. The first class contains methods specifically developed for the consensus problem that operate exclusively with the primal variable $\bbx_{t}$. The general approach here is to perform gradient steps on the local primal variables while also averaging local primal variables with neighbors. The standard first order method is called distributed gradient descent \cite{nedic2009}, and can also be formulated as moving the consensus constraint in \eqref{eq_primal_problem} into the objective function as a penalty term. For better performance in ill-conditioned problems, higher order versions of the primal domain approach include Network Newton \cite{mokhtari2017network} and D-BFGS \cite{eisen2017decentralized}, which employ exact second order and approximated second order information, respectively, to speed up convergence. While all of these methods benefit from achieving at least a linear convergence rate and with low computational cost, they suffer from finding only approximate solutions to \eqref{eq_primal_problem} when using a constant stepsize. They may alternatively use diminishing step sizes to reach the exact solution, but at a slower sublinear convergence rate.

The second class of methods contain those that convert the constrained problem in \eqref{eq_primal_problem} to the dual domain and operate exclusively on the dual variable. These include the standard first order dual descent \cite{cRabbatNowak04} as well as an augmented Lagrangian variation ADMM \cite{Schizas2008-1,Shi2014-ADMM}. Both of these methods perform first order updates in the dual domain and are able to achieve a sub-linear and linear convergence rate, respectively. The D-BFGS method \cite{eisen2017decentralized} performs an approximate second order update in the dual domain. While the convergence rate of these methods is typically not as fast as the linear rate of primal domain methods, they nonetheless improve upon the primal domain methods by finding exact solutions to \eqref{eq_primal_problem}. However, these methods face the additional cost of requiring solutions to internal optimization problems at every iteration of the method. This quality may make them difficult or impractical to use for general objective functions.

The third class of methods combines the faster convergence rate of primal domain methods with the exactness of solutions of the dual domain methods by performing updates on both the primal and dual variables. These may be considered as primal-dual methods. The EXTRA method is an exact first order method with linear convergence rate \cite{shi2015extra} that works effectively as a first order primal-dual method  \cite{mokhtari2016dsa,li2017primal}. The ESOM method \cite{mokhtari2016decentralized} performs a second order update on the primal and variable and first order update on the dual variable to achieve a linear convergence rate to the exact solution without requiring an internal optimization method. The PD-QN method proposed in this work similarly is able to achieve a linear convergence rate without the internal optimization method, but additionally includes an approximate second order update in the dual domain. Thus, the method is able to retain desirable qualities of ESOM without computing exact second order information and providing additional robustness to problems that are ill-conditioned in the dual domain.

\medskip\noindent{\bf Notation remark. } In this paper, we use boldface lower case letters to denote vectors and boldface upper case letters to denote matrices. At any time $t$, the $i$th block of a vector $\bbz_t \in \reals^{np}$ is denoted as $\bbz_{i,t} \in \reals^{p}$, while $\bbz_{n_i.t} \in \reals^{m_i p}$ denotes a concatenation of the components in $i \in n_i$. Likewise,  the $i$th block of matrix $\bbA_t \in \reals^{np \times np}$ is denoted as $\bbA_{i,t} \in \reals^{p \times p}$, while $\bbA_{n_i,t} \in \reals^{m_i p \times m_i p}$ denotes the $(j,k)$ entries of $\bbA$ where $j,k \in n_i$. We further denote by $\bbZ_{\emptyset}$ the matrix containing only the diagonal elements of $\bbZ$.

%
\section{Primal-Dual Quasi-Newton (PD-QN) Method}\label{sec_dbfgs}

We introduce the Primal-Dual Quasi-Newton (PD-QN) algorithm as a fully distributed quasi-Newton update on the augmented Lagrangian function $\ccalL_{\alpha}(\bbx,\bby)$. We refer to it as a primal-dual quasi-Newton method because a quasi-Newton update is used to approximate second order information for both the primal and dual updates. In the primal domain, the second order information is approximated to approximately solve the update in \eqref{eq_mm_p} in a distributed manner, while second order information is approximated in dual domain to make the dual update more robust in ill-conditioned settings and better performing in practice. For notational convenience, we define the primal and dual gradients of $\ccalL_{\alpha}(\bbx_{t},\bby_{t})$, labelled $\bbg_{t}$ and $\bbh_{t}$, respectively, as
\begin{align}\label{eq_primal_grad}
\bbg_{t} &:= \nabla f(\bbx_{t}) + \bby_{t} + \alpha (\bbI-\bbZ)\bbx_{t}, \\
\bbh_{t} &:= (\bbI - \bbZ) \bbx_{t+1}. \label{eq_dual_gradient}
\end{align}
The full gradients $\bbg_t, \bbh_t \in \reals^{np}$ stacks the local gradients at each node, e.g. $\bbg_t = [\bbg_{1,t}; \bbg_{2,t}; \hdots; \bbg_{n,t}]$, where $\bbg_{i,t} \in \reals^p$ is computed locally by node $i$. In the following subsections, we proceed to derive the primal and dual updates of the PD-QN method.

\subsection{Primal update}\label{sec_primal_update}

We seek to replace the non-distributed primal update in \eqref{eq_mm_p} with an update that both allows distributed computation and does not require the explicit solving of a subproblem. To derive such an update, we recall that any descent-based method is in fact the solution of a quadratic Taylor series approximation of the objective function, centered at the current iterate. It is natural then to consider solving the optimization problem in \eqref{eq_mm_p} in the same manner. Consider the quadratic approximation of $\ccalL_{\alpha}(\bbx,\bby_{t})$, centered at $\bbx_{t}$, expressed as
\begin{align} \label{eq_aug_lagrangian_quad}
\hat{\ccalL}_{\alpha}(\bbx,\bby_{t}) &= \ccalL_{\alpha}(\bbx_{t},\bby_{t}) + \nabla \ccalL_{\alpha}(\bbx_{t},\bby_{t})^T(\bbx - \bbx_{t}) \nonumber \\
&\qquad + \frac{1}{2} (\bbx- \bbx_{t}) \bbG_{t}^T(\bbx - \bbx_{t})
\end{align}
where $\bbG_{t}$ is some matrix that approximates the second order information of the augmented Lagrangian. As \eqref{eq_aug_lagrangian_quad} is a quadratic function, the vector $\bbx_{t+1}$ that minimizes it can be found explicitly with a closed form solution. This provides us the following primal variable update
\begin{align}  \label{eq_pdqn_p}
\bbx_{t+1} &= \bbx_{t} - \bbG^{-1}_{t} \bbg_{t}, 
\end{align}
This then provides us the general form of our primal update in PD-QN. A well studied choice of Hessian approximation inverse matrix $\bbG_{t}$ is that given by the quasi-Newton BFGS method \cite{dennis1974characterization,powell1976some}. Consider the particular structure of the Hessian of the augmented Lagrangian $\nabla^2_{\bbx \bbx} \ccalL_{\alpha}(\bbx,\bby_{t}) = \nabla^2 f(\bbx) + \alpha (\bbI-\bbZ)$. The first term, $\nabla^2 f(\bbx)$, is the Hessian of the local objective functions and thus a locally computable block diagonal matrix. The second term, $\alpha (\bbI-\bbZ)$, is not diagonal but has the sparsity pattern of $\ccalG$ and, more importantly, is constant matrix. To construct an approximate matrix then, it is only necessary for nodes to approximate the first term, which can be done solely with local variables and gradients using the standard BFGS quasi-Newton update. To compute this update, define the primal variable and variation variables as
\begin{align} \label{eq_primal_variations}
\bbu_{t} = \bbx_{t+1}-\bbx_{t}, \qquad \bbr_{t} = \bbg_{t+1}-\bbg_{t}.
\end{align}
Because each node $i$ estimates a Hessian that only depends on local variables, its respective BFGS approximation matrix $\bbB_{i,t} \in \reals^{p \times p}$ can be computed using only the local variables and gradients $\bbu_{i,t}$ and $\bbr_{i,t}$. 
Each node maintains and updates at time $t+1$ an approximation of $\nabla^2 f_i (\bbx_i)$ with the iterative BFGS quasi-Newton update
\begin{align} \label{eq_bfgs_p}
\bbB_{i,t+1}
             &= \bbB_{i,t}
               +  \frac{ \bbr_{i,t}  \bbr_{i,t}^{T}}{{\bbu_{i,t}}^{T} \bbr_{i,t}}
                - \frac{\bbB_{i,t} \bbu_{i,t}\bbu_{i,t}^{T}\bbB_{i,t}} 
                       {{\bbu_{i,t}^{T}}\bbB_{i,t}{\bbu_{i,t}}}.
\end{align}
The global objective function Hessian approximation is then defined as the block diagonal matrix combining all local approximations, i.e. $\bbB_{t} := \diag \{ \bbB_{i,t}\}_{i=1}^n \in \reals^{np \times np}$ and the respective full Hessian is approximated as $\bbG_{t} := \bbB_{t} + \alpha(\bbI-\bbZ)$. 

To implement the update in \eqref{eq_pdqn_p} in a distributed manner, the $i$th component of the descent direction $\bbG^{-1}_{t} \bbg_{t}$ must be computable by node $i$ using local information and exchanges with neighbors. More specifically, although $\bbG_{t}$ has the required sparsity pattern of $\ccalG$, its inverse does not. We can, however, approximate the inverse as $\bbG^{-1}_{t,K}$ using $K$ terms of the Taylor series expansion of the inverse $\bbG^{-1}_{t} = (\bbB_{t} + \alpha(\bbI-\bbZ))^{-1}$, written as 
\begin{align} \label{eq_hessian_K}
&\bbG^{-1}_{t,K} :=  \\
& \bbD^{-1/2}_{t} \sum_{k=0}^K \left( \bbD^{-1/2}_{t}  \alpha(\bbI \!\!- \!\!2\bbZ_{\emptyset}\! +\! \bbZ) \bbD_{t}^{-1/2} \right)^k \bbD^{-1/2}_{t}, \nonumber
 \end{align}
where the matrix $\bbD_{t} := \bbB_{t} + 2\alpha(\bbI - \bbZ_{\emptyset})$  contains the block diagonal elements of the approximate Hessian $\bbG_{t}$. The resulting descent update $\bbG^{-1}_{t,K} \bbg_{t}$ can indeed be implemented in a distributed manner at each nodes, with $K+1$ local exchanges needed per iteration to use $K$ terms in the series in \eqref{eq_hessian_K}---see, e.g, \cite{mokhtari2017network}. As larger $K$ will result in a better approximation of the matrix inverse, practical implementations require a tradeoff between accuracy and number of local exchanges required by each node. Each node $i$ can compute its local primal descent component $\bbd_{i,t} := [\bbG^{-1}_{t,K} \bbg_{t}]_i$ and subsequent update using the subroutine displayed in Algorithm \ref{alg_primal}.

%
\begin{algorithm}[t] 
\setstretch{1.35}
{\small\begin{algorithmic}[1]
  \REQUIRE $\{\bbx_{i,\tau},\bbg_{i,\tau}\}_{t,t-1}$,$\bbB_{i,t-1}$, Weights $w_{ij}$ for $j\in n_i$
    \STATE Compute $\bbu_{i,t-1},\bbr_{i,t-1}$, $\bbB_{i,t}$ [cf. \eqref{eq_primal_variations}-\eqref{eq_bfgs_p}]
    \STATE Form $\bbD_i = (1-w_{ii})\bbB_{i,t} + 2(1-w_{ii})\bbI$
     \STATE Initialize $\bbd_{i,t} = -\bbD_i^{-1} \bbg_{i,t}$
     \FOR {$k=0,\hdots,K-1$}
     \STATE Exchange $\bbd_{i,t}$ with neighbors $j \in n_i$
     \STATE  $\bbd_{i,t} = \bbD_i^{-1} [\sum_{j \in n_i} w_{ij} \bbd_{j,t} - \bbg_{i,t}$]
     \ENDFOR
     \STATE Local update
            $\bbx_{i,t+1} = \bbx_{i,t} + \bbd_{i,t}$   
      \RETURN $\bbx_{i,t+1}$, $\bbB_{i,t}$
\end{algorithmic}}
\caption{Primal update for node $i$ at time $t$}
\label{alg_primal}
\end{algorithm}

\subsection{Dual update}\label{sec_dual_update}

We proceed to derive the dual update of the PD-QN method, which replaces the first order dual update in \eqref{eq_dual_ascent_d} with a quasi-Newton update similar to that in \eqref{eq_dual_dbfgs}, i.e.
\begin{equation}\label{eq_pdqn_d}
\bby_{t+1} = \bby_{t} +\alpha\bbH^{-1}_{t} \bbh_{t},
\end{equation}
where $\bbH_{t}$ is the dual Hessian approximation and$\alpha>0$ is the quadratic penalty coefficient in the augmented Lagrangian. As in the primal domain, $\bbH_{t}$ is designed to allow for the distributed computation for dual descent direction $\bbe_{i,t} := [\bbH^{-1}_{t}\bbh_{t}]_i$ to be computed locally at node $i$ using local exchanges. We in particular employ the quasi-Newton dual update used in the dual D-BFGS method \cite{eisen2017decentralized,eisen2016decentralized}, the details of which we discuss here. 

Recall that, in traditional BFGS method, the secant condition $ \bbH_{t+1}  \bbv_{t} =\bbs_{t}$ induces desirable properties onto the approximation matrix $\bbH_{t}$ relating to acceleration and robustness in ill-conditioned problems. We therefore construct a matrix $\bbH_{t}$ that is not only distributable across the network, but maintains the global secant condition. Define the diagonal normalization matrix $\bbUpsilon \in \reals^{np}$ whose $i$th block is $m_i^{-1} \bbI$ and a small scalar regularization parameter $\gamma > 0$. We then define the modified neighborhood variable and gradient variations, $\tbv_{n_i,t} \in \reals^{m_i p}$ and $\tbs_{n_i,t} \in \reals^{m_i p}$, as
\begin{align}
\tbv_{n_i,t} &:= \bbUpsilon_{n_i} \left[ \bby_{n_i,t+1} - \bby_{n_i,t} \right] \label{eq_dbfgs_vars} \\
\tbs_{n_i,t} &:= \bbh_{n_i,t+1} - \bbh_{n_i,t} - \gamma\tbv_{n_i,t}\label{eq_dbfgs_grads}.
\end{align}
The neighborhood variations in \eqref{eq_dbfgs_vars} and \eqref{eq_dbfgs_grads} are modified not just in their locality, but also in the normalization by $ \bbUpsilon_{n_i}$ in \eqref{eq_dbfgs_vars} and regularization by $\gamma\tbv_{n_i,t}$ in \eqref{eq_dbfgs_grads}. As $\tbv_{n_i,t}$ and $\tbs_{n_i,t}$ can be obtained locally at each node $i$ with one hope excahnges, each node $i$ computes and maintains a local Hessian approximation $\bbC_{n_i,t} \in \reals^{m_ip \times m_ip}$, which is updated as the solution of a local optimization problem,
\begin{alignat}{2}\label{eq_dbfgs_update}
\bbC_{n_i,t+1} := &\argmin_{\bbZ}\ &&
\text{tr}[ (\bbC_{n_i,t})^{-1} (\bbZ - \gamma \bbI)] - \\\nonumber &&&\qquad\quad
 \text{logdet}[(\bbC_{n_i,t})^{-1} (\bbZ - \gamma \bbI)] - n\\\nonumber
&  \text{   s.t.}\quad &&
\bbZ \tbv_{n_i,t} = \bbs_{n_i,t},\quad \bbZ \succeq \bb0.
\end{alignat}
To update in \eqref{eq_dbfgs_update} provides an updated approximation matrix $\bbC_{n_i,t+1}$ that has eigenvalues greater than $\gamma$ while satisfying a \emph{modified} local secant condition with respect to the normalized variable variation $\tbv_{n_i,t}$. A closed form solution to \eqref{eq_dbfgs_update} exists \cite[Proposition 1]{mokhtari2014res}, and is given by
\begin{align}\label{eq_bfgs_dist}
\bbC_{n_i}(t + 1) &=  \bbC_{n_i,t}+ \frac{\tbs_{n_i,t} \tbs^T_{n_i,t}}{\tbs^T_{n_i,t} \tbv_{n_i,t}}   \\
& \qquad - \frac{\bbC_{n_i,t} \tbv_{n_i,t} \tbv^T_{n_i,t} \bbC_{n_i,t}}{\tbv^T_{n_i,t} \bbC_{n_i,t} \tbv_{n_i,t}} + \gamma \bbI.  \nonumber
\end{align}
Note that the dual BFGS update differs from the more traditional primal BFGS update in \eqref{eq_bfgs_p} both in the use of modified neighborhood variations $\tbv_{n_i,t} $ and $\tbs_{n_i,t}$ and the addition of regularization parameter $\gamma \bbI$. 

\begin{figure}\centering

\def \thisplotscale {0.6}
\def \unit {\thisplotscale cm}

\tikzstyle {block}        = [draw, very thin,
                             rectangle, 
                             minimum height = 1*\unit,
                             minimum width  = 1*\unit]                       

\tikzstyle {blue block}   = [block,
                             fill = blue!20]

\tikzstyle {green  block} = [block,
                             fill = green!20]

\tikzstyle {red block}    = [block,
                             fill = red!20]

\def \deltalabel{ 0.5}

{\fontsize{7}{7}\selectfont\begin{tikzpicture}[x = 1*\unit, y=1*\unit, 
                          shorten >=2pt, shorten <=2pt]

%
%
\path (0,0)        node [red block       ] (Bine) {};
\path (Bine.south) node [red block, below] (Bise) {};
\path (Bise.west)  node [red block, left ] (Bisw) {};
\path (Bisw.north) node [red block, above] (Binw) {};
\path (Bine.north west) ++ (0, \deltalabel) node {$(\bbC^i)^{-1}$}; 
%
%
\path (Binw.west) ++ (-1,0)   node [red   block, left] (rnin) {$\bbh_i$};
\path (Bisw.west) ++ (-1,0)   node [blue  block, left] (rnis) {$\bbh_j$};
\path (rnin.west) ++ (-0.3,0) node [red   block, left] (vnin) {$\bby_i$};
\path (rnis.west) ++ (-0.3,0) node [blue  block, left] (vnis) {$\bby_j$};
\path (rnin.north) ++ (0, \deltalabel) node {$\bbh_{n_i}$}; 
\path (vnin.north) ++ (0, \deltalabel) node {$\bby_{n_i}$}; 
%
%
\path (Bine.east) ++ (1,0) node [red   block, right] (gnin) {$\bbh_i$};
\path (Bise.east) ++ (1,0) node [blue  block, right] (gnis) {$\bbh_j$};
\path (gnin.north) ++ (0, \deltalabel) node {$\bbh_{n_i}$}; 
%
%
\path (gnin.east) ++ (1,0) node [red   block, right] (dnin) {$\bbe_i^i$};
\path (gnis.east) ++ (1,0) node [blue  block, right] (dnis) {$\bbe_j^i$};
\path (dnin.north) ++ (0, \deltalabel) node {$\bbe^i_{n_i}$}; 
%
%
\path (dnin.south east) ++ (1,0) node [red block, right] (dil) {$\bbe_i^i$};
\path (dil.east)        ++ (1,0) node [red block, right] (dir) {$\bbe_i^j$};
\path (rnin.south) -- (Binw.south)  node [midway] {$\text{\eqref{eq_bfgs_dist}} \atop
                                                    \Rightarrow$};
\path (Bine.south) -- (gnin.south)  node [midway] {$\times$};
\path (gnin.south) -- (dnin.south)  node [midway] {$=$};
\path (dil.east)   -- (dir.west)    node [midway, draw, circle, inner sep=1] (plus) {$+$};
\path[draw, -stealth] (plus) -- ++(0,1.5) node [red block, above] {$\bbe_i$};;

%
%
%
\path (0, -3.3)    node [blue block       ] (Bjne) {};
\path (Bjne.south) node [blue block, below] (Bjse) {};
\path (Bjse.west)  node [blue block, left ] (Bjsw) {};
\path (Bjsw.north) node [blue block, above] (Bjnw) {};
\path (Bjse.south west) ++ (0,-\deltalabel) node {$(\bbC^j)^{-1}$}; 
%
%
\path (Bjnw.west) ++ (-1,0)   node [red   block, left] (rnjn) {$\bbh_i$};
\path (Bjsw.west) ++ (-1,0)   node [blue  block, left] (rnjs) {$\bbh_j$};
\path (rnjn.west) ++ (-0.3,0) node [red   block, left] (vnjn) {$\bby_i$};
\path (rnjs.west) ++ (-0.3,0) node [blue  block, left] (vnjs) {$\bby_j$};
\path (rnjs.south) ++ (0,-\deltalabel) node {$\bbh_{n_j}$}; 
\path (vnjs.south) ++ (0,-\deltalabel) node {$\bby_{n_j}$}; 
%
%
\path (Bjne.east) ++ (1,0) node [red   block, right] (gnjn) {$\bbh_i$};
\path (Bjse.east) ++ (1,0) node [blue  block, right] (gnjs) {$\bbh_j$};
\path (gnjs.south) ++ (0,-\deltalabel) node {$\bbh_{n_j}$}; 
%
%
\path (gnjn.east) ++ (1,0) node [red   block, right] (dnjn) {$\bbe_i^j$};
\path (gnjs.east) ++ (1,0) node [blue  block, right] (dnjs) {$\bbe_j^j$};
\path (dnjs.south) ++ (0,-\deltalabel) node {$\bbe^j_{n_j}$}; 
%
%
\path (dnjn.south east) ++ (1,0) node [blue block, right] (djl) {$\bbe_j^j$};
\path (djl.east)        ++ (1,0) node [blue block, right] (djr) {$\bbe_j^i$};
%
%
\path (rnjn.south) -- (Bjnw.south)  node [midway] {$\text{\eqref{eq_bfgs_dist}} \atop
                                                    \Rightarrow$};
\path (Bjne.south) -- (gnjn.south)  node [midway] {$\times$};
\path (gnjn.south) -- (dnjn.south)  node [midway] {$=$};
\path (djl.east)   -- (djr.west)    node [midway, draw, circle, inner sep=1] (plus) {$+$};
\path[draw, -stealth] (plus) -- ++(0,-1.5) node [blue block, below] {$\bbe_j$};;

%
%
\path (dnin.east) edge [-stealth, black!99, bend left]  (dil.north);
\path (dnis.east) edge [-stealth, black!99, bend left]  (djr.north);
\path (dnjs.east) edge [-stealth, black!99, bend right] (djl.south);
\path (dnjn.east) edge [-stealth, black!99, bend right] (dir.south);

%
%
%
\path (rnin.east) edge [-stealth, black!99, bend left] (rnjn.east);
\path (rnjs.west) edge [-stealth, black!99, bend left] (rnis.west);
\path (vnin.east) edge [-stealth, black!99, bend left] (vnjn.east);
\path (vnjs.west) edge [-stealth, black!99, bend left] (vnis.west);

\end{tikzpicture}}
\caption{PD-QN dual variable flow. Nodes exchange variable and gradients -- $\bby_i$ and $\bbh_i$ sent to $j$ and $\bby_j$ and $\bbh_j$ sent to $i$ -- to build variable and gradient variations $\tbv$ and $\tbs$ that they use to determine local curvature matrices -- $\bbC_{n_i}$ and $\bbC_{n_j}$. They then use gradients $\bbh_{n_i}$ and $\bbh_{n_j}$ to compute descent directions $\bbe_{n_i}^i$ and  $\bbe_{n_j}^j$. These contain a piece to add locally -- $\bbe_i^i$ stays at node $i$ and $\bbe_j^j$ stays at node -- and a piece to add at neighbors -- $\bbe_j^i$ is sent to node $j$ and $\bbe_i^j$ is sent to node $i$.}
\label{fig_variable_flow_diagram} \end{figure}
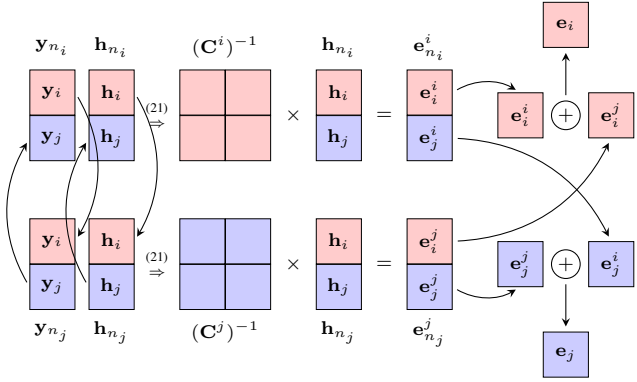

Observe that, through the use of \emph{neighborhood} variables in the update in \eqref{eq_bfgs_dist}, nodes approximate the dual Hessian of themselves \emph{and} their neighbors. They subsequently use $\bbC_{n_i,t}$ along with an additional small regularization parameter $1 \geq \Gamma >0$ to compute the neighborhood descent direction $\bbe^i_{n_i,t} \in \reals^{m_i p}$ as
\begin{equation}
\bbe^i_{n_i,t} = - \left( \bbC_{n_i,t}^{-1} + \Gamma \bbUpsilon_{n_i} \right) \bbh_{n_i,t}.
\label{eq_direction_local}
\end{equation}
The neighborhood descent direction $\bbe^i_{n_i,t} \in \reals^{m_ip}$ contains components for variables of node $i$ itself and all neighbors $j \in n_i$ -- see Fig. \ref{fig_variable_flow_diagram}. Likewise, neighboring nodes $j \in n_i$ contain a descent component of the form $\bbe^j_{i,t}$. The local descent $\bbe_{i,t}$ is then given by the sum of the components $\bbe^j_{i,t}$ for all neighbors $j\in n_i$, i.e. $\bbe_{i,t} = \sum_{j \in n_i}  \bbe^j_{i,t}$.

The global dual Hessian approximation $\bbH_{t}$ in \eqref{eq_pdqn_d} can be derived from all nodes performing the local update in \eqref{eq_direction_local} simultaneously in parallel. More precisely, it can be shown that $\bbH^{-1}_{t} = \hbH^{-1}_{t} + \Gamma\bbI$, where $\hbH$ satisfies the global secant condition $ \hbH_{t+1}  \bbv_{t} =\bbs_{t}$---see \cite[Proposition 1]{eisen2017decentralized} for details. Moreover, the update in \eqref{eq_pdqn_d} can be computed distributedly using a sequence of local exchanges. The necessary exchanges for node $i$ are detailed in Algorithm \ref{alg_dual}.

%
\begin{algorithm}[t] 
\setstretch{1.35}
\begin{algorithmic}[1]\small
  \REQUIRE $\{\bby_{n_i,\tau},\bbh_{n_i,\tau}\}_{t,t-1}$,$\bbC_{n_i,t-1}$, $\eps_d$,$\gamma,\Gamma >0$
        \STATE Compute $\tbv_{ n_i,t-1},\tbs_{ n_i,t-1},\bbC_{n_i,t}$ [cf.\eqref{eq_dbfgs_vars}--\eqref{eq_dbfgs_update}] 
  	\STATE Compute $\bbe^i_{n_i,t} = - (\bbC_{n_i,t}^{-1} +\Gamma\bbUpsilon_{n_i}) \bbh_{n_i,t}$ [cf. \eqref{eq_direction_local}]
  	\STATE Exchange $\bbe^i_{j,t}$ with neighbors $j \in n_i$
  	\STATE Compute descent dir. $\displaystyle{\bbe_{i,t} := \sum_{j \in  n_i}  \bbe^j_{i,t}}$
  	\STATE Update $\bby_{i,t+1} = \bby_{i,t} +\eps_d\bbe_{i,t}$
      \RETURN $\bbx_{i,t+1}$, $\bbC_{n_i,t}$
\end{algorithmic}
\caption{Dual update for node $i$ at time $t$}
\label{alg_dual}
\end{algorithm}

%

%
For the ease of presentation, a full description of variable exchanges necessary to perform the update is not shown in  Algorithm \ref{alg_dual}. We present in Fig. \ref{fig_variable_flow_diagram} a diagram of the flow of variables among neighbors. Variable and gradients are exchanged -- $\bby_{i,t}$ and $\bbh_{i,t}$ are sent to node $j$ and $\bby_{j,t}$ and $\bbh_{j,t}$ are sent to node $i$ -- and \eqref{eq_bfgs_dist} is used to compute the curvature estimation matrices $\bbC_{n_i,t}$ and $\bbC^j_{t}$. Using these Hessian approximations, the nodes premultiply the inverse by the neighborhood gradients $\bbh_{n_{i,t}}$ and $\bbh_{n_{j,t}}$ to obtain the neighborhood descent directions -- $\bbe^i_{n_i,t}$ and  $\bbe^j_{n_j,t}$. These descent directions contain a piece to be added locally -- $\bbe^i_{i,t}$ stays at node $i$ and $\bbe^j_{j,t}$ stays at node -- and a piece to be added at the neighboring node -- $\bbe^i_{j,t}$ is sent to node $j$ and $\bbe^j_{i,t}$ is sent to node $i$. The local descent direction $\bbe_{i,t}$ is the addition of the locally computed $\bbe^i_{i,t}$ and the remotely computed $\bbe^j_{i,t}$.

We conclude the discussion on the dual update with a brief series of remarks.

\begin{remark}\label{remark_comm}\normalfont

As can be seen in the variable flow of Fig. \ref{fig_variable_flow_diagram}, the dual update of PD-QN requires a larger amount of information exchanges than is needed for other consensus methods. Indeed, the total communication overhead for PD-QN includes $K+1$ exchanges for the primal update and 4 exchanges for the dual update---$K+5$ in total. This is in comparison to the 2 exchanges needed for first order methods like DA and ADMM and $K+3$ exchanges needed for ESOM. While this burden is higher, the gains made by PD-QN in iterations necessary to convergence can, in many cases, still render PD-QN a preferable approach to alternatives. This tradeoff is explored in the numerical results in Section \ref{sec_numerical_results} of this paper.
\end{remark}

%
\begin{remark}\label{rmk_inner_product_negative}\normalfont 
For the problem in \eqref{eq_dbfgs_update} to have a solution and the update in \eqref{eq_bfgs_dist} to be valid the inner product between the neighborhood variations must be $\tbv_{n_i,t}^T\tbs_{n_i,t} > 0$. This condition imposes a restriction in functions that can be handled by PD-QN. In practical implementations, however, we can check the value of this inner product and proceed to update $\bbC_{n_i,t}$ only when it satisfies $\tbv_{n_i,t}^T\tbs_{n_i,t} > 0$.  \end{remark}

%

\subsection{PD-QN Summary}

Here we summarize the full details of the PD-QN update, combining both the primal and dual updates in \eqref{eq_pdqn_p} and \eqref{eq_pdqn_d} respectively, restated here as 
\begin{align}
\bbx_{t+1} &= \bbx_{t} -\bbG^{-1}_{t,K} \bbg_{t}, \label{eq_pdqn_final1} \\
\bby_{t+1} &= \bby_{t} + \alpha \bbH^{-1}_{t} \bbh_{t}, \label{eq_pdqn_final2} 
\end{align}
A summary of the variables relevant to these updates are provided in Table \ref{tab_variables} for reference.

\begin{table}[t]
\centering
\caption{PD-QN variable summary}
\label{tab_variables}
\begin{tabular}{l | ll}
                        & Primal   & Dual     \\ \hline
Variable                & $\bbx$   & $\bby$  \\
Gradient                & $\bbg$   & $\bbh$   \\
Variable variation      & $\bbu$   & $\bbv$   \\
Gradient variaton       & $\bbr$   & $\bbs$   \\
Local Hessian approx.   & $\bbB_i$ & $\bbC_{n_i}$ \\
Global Hessian approx.  & $\bbG$   & $\bbH$   \\
Local Descent direction & $\bbd_i$ & $\bbe_i$
\end{tabular}
\end{table}

The complete PD-QN algorithm at node $i$, including both the primal and dual updates, is summarized in Algorithm \ref{alg_pdqn}.  Each node begins with initial variables $\bbx_{i,0}$ and $\bby_{i,0}$ and exchanges variables with neighbors. For each step $t$, the node computes its primal gradient in Step 3 and update the primal variable in Step 4 using Algorithm \ref{alg_primal}. They exchange the updated variables with neighbors in Step 5 and use them to compute and exchange dual gradients in Steps 6 and 7. The dual variable is updated with Algorithm \ref{alg_dual} in Step 8, after which they are exchanged in Step 9.

\section{Convergence Analysis}\label{sec_convergence}

We analyze the convergence of PD-QN method performed on the consensus optimization problem in \eqref{eq_primal_problem}. To begin, we make the following assumption on the eigenvalues of the objective function Hessian,
 \begin{assumption} \label{as_strongly_convex}
The aggregate objective function $f(\bbx) = \sum_{i=1}^n f_i(\bbx_i)$ is twice differentiable and the eigenvalues of the objective function Hessian are nonnegative and bounded from above and below by positive constants $0 < \mu < L < \infty $, i.e. 
\begin{equation}
\mu \bbI \preceq \nabla^2 f(\bbx) \preceq L\bbI.
\label{eq_hessian_bounds_sc}
\end{equation}
\end{assumption} 

 The upper bound $L$ on the eigenvalues of the Hessian in Assumption \ref{as_strongly_convex} implies that the associated gradient $\bbg(\bbx)$ is Lipschitz continuous with parameter $L$, i.e. $\| \bbg(\bbx) - \bbg(\bbx ')\| \leq  L\|\bbx - \bbx'\|$. The aggregate objective functions is additionally strongly convex with parameter $\mu$. These are often considered standard assumptions in distributed convex optimization methods. Note that, the weak convexity (i.e. $\mu=0$) is sometimes analyzed as well, but is not considered in this paper. Many problems in distributed machine learning with weakly convex objective functions are often supplemented with a strongly convex regularizer (e.g. $\ell$-2 norm). We further make an assumption on the eigenvalues of the primal Hessian approximation matrix $\bbB_{t}$.
 
 \begin{assumption} \label{as_bfgs_bound}
There exist positive constants $0 < \psi < L < \Psi $ such that the eigenvalues of the primal Hessian approximation matrix $\bbB_{t}$ are bounded from above and below as
\begin{equation}
\psi \bbI \preceq \bbB_{t} \preceq \Psi\bbI.
\label{eq_bfgs_bound}
\end{equation}
\end{assumption} 
\begin{remark}\normalfont
The bounds imposed on the eigenvalues of $\bbB_{t}$ are, in general, not standard assumptions. While the matrix $\bbB_{t}$ is guaranteed to be positive definite, the lower eigenvalue can be arbitrarily small. However, there are a number of techniques commonly used to satisfy this assumption in practice.  These include both adding small regularization terms to both $\bbB_{t}$ and $\bbB^{-1}_{t}$---see, e.g. \cite{mokhtari2014res}---and using the popular limited memory version of the BFGS update in \eqref{eq_bfgs_p}, which induces the necessary bounds in \eqref{eq_bfgs_bound}---see, e.g., \cite{mokhtari2015global}. In this paper, we assume such bounds exist for the ease of analysis. We also observe in the numerical experiments in Section \ref{sec_numerical_results} that such regularization techniques are often not necessary in practice.
\end{remark}

%
\begin{algorithm}[t!] 
\setstretch{1.35}
\small{\begin{algorithmic}[1]
  \REQUIRE
  $\bbx_{i,0}, \bby_{i,0}, \bbB_{i,0}, \bbC_{n_i,0}, \eps_d, \alpha$
  \STATE Exchange initial variables with neighbors $j \in n_i$
  \FOR{$t = 0,1,2, \hdots$}     
     \STATE Grad. $\bbg_{i,t} = \bbx_{i,t} - \sum_{j \in n_i} \!\!\!w_{ij}\bbx_{j,t} + \bby_{i,t} + \alpha \nabla f_i(\bbx_{t})$
     \STATE Update primal $\bbx_{i,t+1}$ with Algorithm \ref{alg_primal}.
     \STATE Exchange $\bbx_{i,t+1}$ with neighbors $j \in n_i$
       \STATE Grad. $\bbh_{i,t} = \bbx_{i,t+1} - \sum_{j \in n_i} w_{ij} \bbx_{j,t+1}$
       \STATE Exchange $\bbh_{i,t}$ with neighbors $j \in n_i$
       \STATE Update dual $\bby_{i,t+1}$ with Algorithm \ref{alg_dual}
       \STATE Exchange $\bby_{i,t+1}$ with neighbors $j \in n_i$
  \ENDFOR
\end{algorithmic}}
\caption{PD-QN method at node $i$}
\label{alg_pdqn}
\end{algorithm}
 
 \begin{assumption} \label{as_inner_product}
For all $i$ and $t$, the inner product between the neighborhood dual variable and gradient vector variations is strictly positive, i.e. $\tbv_{n_i,t}^T \tbs_{n_i,t} > 0$.
\end{assumption} 
We state this assumption explicitly to ensure all local dual Hessian approximations $\bbC_{n_i,t}$ are well defined in \eqref{eq_bfgs_dist}. While this may not hold in practice, we may set $\bbC_{n_i,t+1} = \bbC_{n_i,t}$---see Remark \ref{rmk_inner_product_negative}---which we stress does not have any bearing on the proceeding analysis.

Before deriving the primary theoretical results of the PD-QN method, we first establish some properties of the Hessian approximation matrices for the primal and dual domains, denoted $\bbG_{t,K}$ and $\bbH_{t}$. The following lemmata characterize the eigenvalues of the $K$-th order inverse approximation of the primal Hessian approximation $\bbG^{-1}_{t,K}$ and the dual Hessian inverse approximation $\bbH^{-1}_{t}$, respectively.

\begin{lemma}\label{lemma_G_bound}
Consider the primal update in the PD-QN update introduced in \eqref{eq_pdqn_p}-\eqref{eq_hessian_K}. If Assumptions \ref{as_weight_bound}-\ref{as_bfgs_bound} hold, then the eigenvalues of the primal Hessian inverse approximation $\bbG^{-1}_{t,K}$ are uniformly bounded as 
\begin{equation}
\lambda \bbI \preceq \bbG^{-1}_{t,K} \preceq  \Lambda \bbI,
\label{eq_prop_eigen_bounds_primal} 
\end{equation} 
where the constants $\lambda$ and $\Lambda$ are defined as
\begin{align}
\lambda := \frac{1}{2\alpha(1-\delta) + \Psi}, \quad \Lambda := \frac{1-\rho^{K+1}}{(1-\rho)(2\alpha(1-\Delta)+\psi)} \nonumber
\end{align}
and $ \rho:= (2\alpha(1-\delta))/(2\alpha(1-\delta)+\psi)$.
\end{lemma}

\begin{myproof}
The proof can be found for the similar result in \cite[Lemma 2]{mokhtari2017network} and is excluded here for space considerations.
 \end{myproof}
\begin{lemma}\label{lemma_H_bound}
Consider the dual update in the PD-QN method introduced in \eqref{eq_pdqn_d}-\eqref{eq_direction_local}. Further, recall both the positive constants $\gamma$ and $\Gamma \leq 1$ as the regularization parameters of dual Hessian and the definition of its global approximation $\bbH^{-1}_{t} =\hbH^{-1}_{t} + \Gamma \bbI$. If Assumption \ref{as_inner_product} holds, the eigenvalues of the dual Hessian inverse approximation $\bbH^{-1}_{t}$ are bounded as 
\begin{equation}
\Gamma \bbI \preceq \bbH^{-1}_{t} \preceq  P \bbI,
\label{eq_prop_eigen_bounds_dual} 
\end{equation} 
where $P := \left( \Gamma + n/\gamma \right)$ and $n$ is the size of network. 
\end{lemma}
 
\begin{myproof}
To establish the lower bound in \eqref{eq_prop_eigen_bounds_dual}, consider that $\hbH^{-1}_{t}$ is a sum of positive semidefinite matrices and is therefore a positive semidefinite matrix with eigenvalues greater than or equal to 0. The upper bound, on the other hand, follows from the fact that each $\hbH_{t}$ is the sum of $n$ matrices, where $\bbC_{n_i,t}^{-1} \preceq 1/\gamma \bbI$ for all $i$. Adding the regularization term $\Gamma\bbI$ provides the upper bound in\eqref{eq_prop_eigen_bounds_dual}.
 \end{myproof} 
 
In Lemmata \ref{lemma_G_bound} and \ref{lemma_H_bound} we show that there exists lower and upper bounds on the eigenvalues of both the primal and dual Hessian inverse approximation matrix. From here, we proceed to demonstrate the linear convergence of the PD-QN method. The following lemma establishes an important relationship between the primal and dual variables using the PD-QN updates.

\begin{lemma}\label{lemma_pdqn_error}
Consider the updates of PD-QN in \eqref{eq_pdqn_final1} and \eqref{eq_pdqn_final2}, where we recall the approximate primal and dual Hessian inverses $\bbG^{-1}_{t,K}$ and $\bbH^{-1}_{t}$. If Assumptions \ref{as_weight_bound}-\ref{as_inner_product} hold, then primal and dual iterates generated by PD-QN satisfy
\begin{align}\label{opt_res_PD-QN2}
&\nabla f(\bbx_{t+1})-\nabla f(\bbx^*) + \bby_{t+1}-\bby^* + \bbsigma_{t}=\bb0,
\end{align}
where the error vector $\bbsigma_{t}$ is defined as 
\begin{align}\label{esom_error_vec}
\bbsigma_{t} &:= \nabla f(\bbx_{t})-\nabla f(\bbx_{t+1}) - \alpha\bbH_{t}^{-1}(\bbI-\bbZ)(\bbx_{t+1}-\bbx^*) 
 \nonumber\\ &\qquad
 + \left[\bbG_{t,K}-\alpha(\bbI-\bbZ)\right](\bbx_{t+1}-\bbx_{t}).
\end{align}
\end{lemma}
\begin{myproof}
See Appendix \ref{app_lemma_pdqn_error}.
 \end{myproof}

In Lemma \ref{lemma_pdqn_error}, we establish a relationship between the primal and dual variables, that is similar to one used in the convergence of the method of multipliers---see \cite{mokhtari2016decentralized}. The PD-QN method includes an additional error term $\bbsigma_{t}$ that encompasses two modifications to MM: (i) the use of the approximate primal Hessian approximation $\bbG_{t}$ rather than the true primal Lagrangian Hessian $\nabla^2 \ccalL(\bbx,\bbnu_{t})$ and (ii) the use of the dual quasi-Newton matrix $\bbH_{t}$ rather than a first order dual update.

From here, we may establish a convergence rate of the PD-QN method, first by establishing a linear convergence rate of a Lyapunov function. To define the Lyapunov function, we first define an appended variable $\bbz_t$ and matrix $\bbJ_t$ as
\begin{equation}\label{eq_aug_defs}
\bbz_t = \begin{bmatrix} 
\bbx_t  \\
\bbnu_t
\end{bmatrix} \quad
\bbJ_t = \begin{bmatrix} 
\alpha \bbR_t  & \bb0  \\
\bb0 & \bbH_t \end{bmatrix}. 
\end{equation}


\begin{theorem}\label{thm:esom_linear_convg}
Consider PD-QN as introduced in \eqref{eq_pdqn_final1}-\eqref{eq_pdqn_final2}. Consider arbitrary constants $\beta>1$ and $\phi>1$ and $\zeta$ as a positive constant. Further, recall the definitions of the vector $\bbz_t$ and matrix $\bbJ_t$ in \eqref{eq_aug_defs} and consider $\hat{\delta}$ as the smallest non-zero eigenvalue of the matrix $\bbI-\bbZ$. If Assumptions \ref{as_weight_bound}-\ref{as_inner_product} hold, then the sequence of Lyapunov functions $\|\bbz_{t}-\bbz^*\|_{\bbJ._t}$ generated by PD-QN satisfies 
\begin{equation}\label{pdqn_lin_convg}
\|\bbz_{t+1}-\bbz^*\|_{\bbJ_t}^2 \ \leq\  \frac{1}{1+\kappa_t} \ \|\bbz_{t}-\bbz^*\|_{\bbJ_t}^2.
\end{equation}
where the sequence $\kappa_t$ is given by
\begin{align}
&\kappa_t=\min \Bigg\{ 	    
        \left( \frac{\beta^2}{P(\beta-1)\hat{\delta}} - \frac{2 \beta \phi \Gamma^2}{P (\phi-1)\hat{\delta}}\right)^{-1}\left(\alpha\Sigma - 2\alpha\zeta L^2/\Sigma\right), \nonumber
	\\
	&
	\frac{2\alpha\hat{\delta}}{\phi\beta(\mu+L)},
	\left( \Sigma - \frac{2 \beta \phi \alpha}{P(\phi-1)\hat{\delta}} \right)^{-1}\!\! \left(\frac{2\mu L}{\mu + L}\! -\! \frac{1}{\zeta}\! -\! \frac{4\alpha^2 P \zeta}{(1-\delta)^{-1}} \right)
	\Bigg\}.\nonumber
\end{align}
and $\Sigma := \Lambda^{-1} - 2\alpha(1-\delta)$.
\end{theorem}
\begin{myproof}
See Appendix \ref{app_thm_esom_linear_convg}.
 \end{myproof}

 Theorem \ref{thm:esom_linear_convg} provides a linear convergence rate of the PD-QN method in terms of the sequence $\| \bbz_t - \bbz^*\|_{\bbJ_2}^2$. From here, it remains to show that the sequence of primal $\bbx_t$ also converges to the optimal argument $\bbx^*$ at a linear rate. The resulting corollary follows as a direct consequence of the preceding theorem, and is presented below.

\begin{corollary}\label{esom_approx_error2}
If Assumptions \ref{as_weight_bound}-\ref{as_inner_product} hold, then the sequence of squared errors $\|\bbx_{t}-\bbx^*\|^2$ generated by PD-QN converges to zero at a linear rate, i.e.,
\begin{equation}\label{ESOM_lin_convg_22}
\|\bbx_{t}-\bbx^*\|^2  \leq  \left(\frac{1}{1+\min_{t}\{\kappa_t\}}\right)^t  \frac{\|\bbz_0-\bbz^*\|_{\bbJ_t}^2}{\alpha \Sigma}.
\end{equation}
\end{corollary}
\begin{myproof}
According to the definition of the sequence $\bbz_t$ and matrix $\bbJ_2$, we can write $\|\bbz_t-\bbz^*\|_{\bbJ_t}^2=\alpha\|\bbx_t-\bbx^*\|_{\bbR_t}^2+ \|\bbnu_t-\bbnu^*\|_{\bbH_t}^2$ and subsequently lower bounded as $\|\bbz_t-\bbz^*\|_{\bbJ_t}^2 \geq \alpha \Sigma \|\bbx_t-\bbx^*\|^2+ (1/P) \|\bbnu_t-\bbnu^*\|^2$ which implies that $\|\bbx_t-\bbx^*\|^2\leq (1/\alpha \Sigma)\|\bbz_t-\bbz^*\|_{\bbJ_t}^2$. Considering this result and linear convergence of the sequence $\|\bbz_t-\bbz^*\|_{\bbJ_t}^2$ in \eqref{pdqn_lin_convg}, the claim in \eqref{ESOM_lin_convg_22} follows. 
\end{myproof}

 The results here establish a linear convergence rate to the exact solution of the consensus problem---this rate is comparable with state of the art methods such as EXTRA \cite{shi2015extra} and ESOM \cite{mokhtari2016decentralized}. Furthermore, we stress that PD-QN does not require internal minimization steps used in pure dual methods or computation of Hessian information used in pure second order methods. We proceed to show the performance of PD-QN in numerical studies.

%
\section{Numerical Results} \label{sec_numerical_results}

\begin{figure*}[t]
\centering
	\begin{subfigure}[t]{.45\textwidth}
		\centering
		\includegraphics[height=\textheight,width=\textwidth,keepaspectratio]{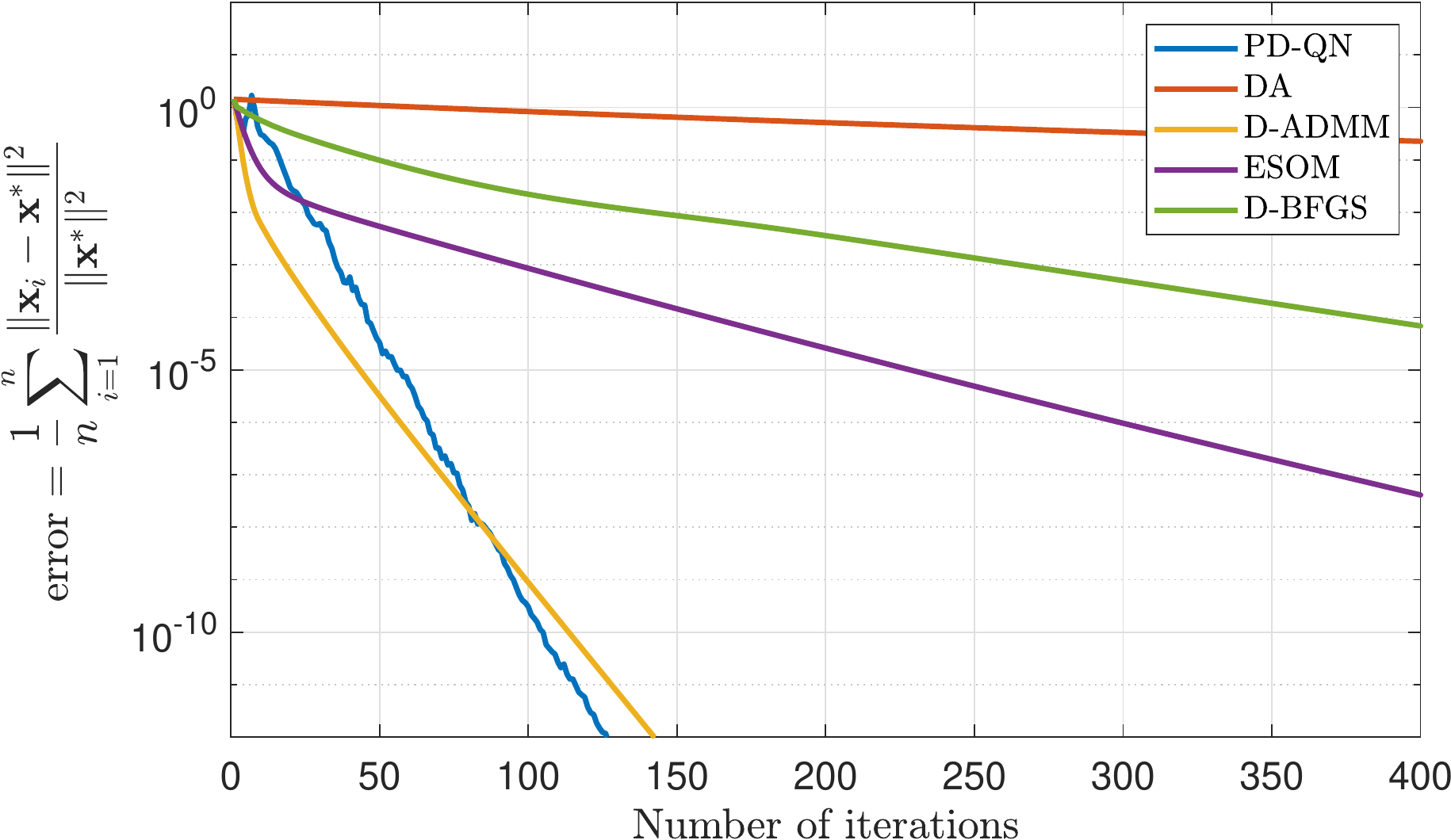}
		\caption{}\label{fig:4a}
	\end{subfigure}
	\begin{subfigure}[t]{.45\textwidth}
		\centering
		\includegraphics[height=\textheight,width=\textwidth,keepaspectratio]{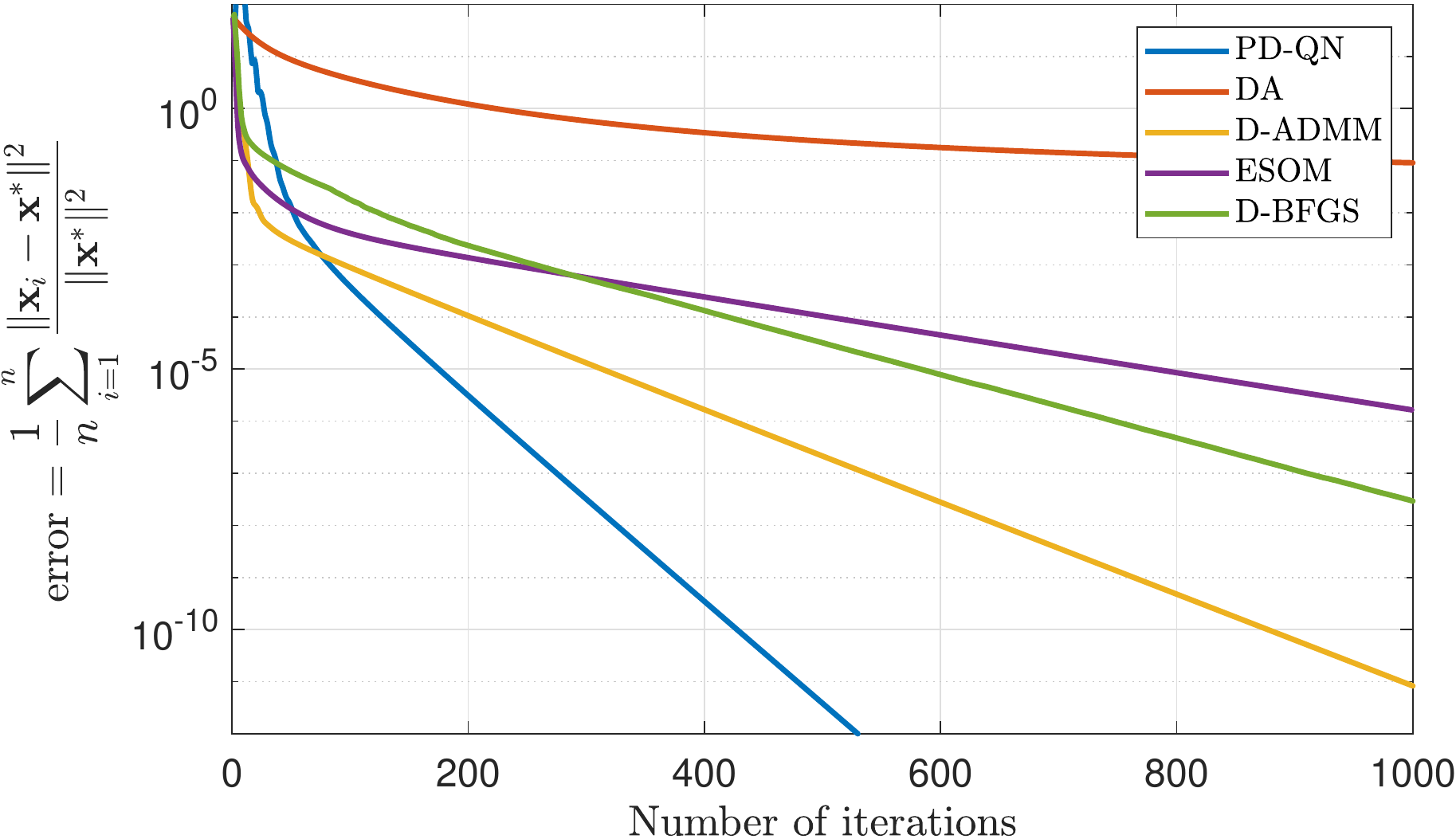}
		\caption{}\label{fig:4b}
	\end{subfigure} 	
	\caption{Convergence paths for exact distributed methods (a) small and (b) large condition number. PD-QN provides significant improvement in convergence time over other methods in both cases.}\label{figure_simulation_quad}
\end{figure*}

We provide numerical simulations of the performance of PD-QN on the consensus problem in \eqref{eq_primal_problem} to compare to other first and second order distributed methods that solve for the exact solution. In particular, we compare against first order methods dual ascent (DA) \cite{cRabbatNowak04} and D-ADMM \cite{Schizas2008-1,Shi2014-ADMM}, and second order methods ESOM \cite{mokhtari2016decentralized} and D-BFGS \cite{eisen2017decentralized}. We demonstrate these results for two common objective functions in distributed learning---linear least squares regression and logistic regression. We begin with the linear least squares regression problem, which is well known to be formulated as a quadratic program. Specifically, we consider the following objective function
\begin{align}
f(\bbx) := \sum_{i=1}^n \frac{1}{2} \bbx^T \bbA_i \bbx + \bbb_i^T \bbx
\label{eq_simulation_problem}
\end{align}
where $\bbA_i \in \reals^{p \times p}$ and $\bbb_i \in \reals^p$ are parameters available to node $i$.  As a means of controlling the condition number of the problem, we define the matrices $\bbA_i = \text{diag}\{ \bba_i \}$, and for a chosen condition number $10^{2\eta}$, $\bba_i$ is given $p/2$ elements chosen randomly from the interval $[1, 10^1, \hdots, 10^{\eta}]$ and $p/2$ elements chosen randomly from the interval $[1, 10^{-1}, \hdots, 10^{-\eta}]$. It is then the case that the full matrix $\sum_{i=1}^n \bbA_i$ has eigenvalues in the range $[n 10^{-\eta}, n 10^{\eta}]$ and the intended condition number. The vectors $\bbb_i$, alternatively, are chosen uniformly and randomly from the box $[0,1]^p$.

In all initial simulations we fix the variable dimension $p=5$ and $n=20$ nodes and use a $d$-regular cycle for the graph, in which $d$ is an even number and nodes are connected to their $d/2$ nearest neighbors in either direction. The others parameters such as condition number $10^{\eta}$ and and number of nodes $n$ are varied by simulation. The regularization parameters for BFGS are chosen to be $\gamma = \Gamma = 10^{-1}$. For all methods we choose a constant stepsize and attempt to pick the largest stepsize for which the algorithms are observed to converge. Representative convergence paths for the five compared methods are shown in Figure \ref{figure_simulation_quad}. Specifically, the convergence paths represent the relative error with respect to the exact solution $\bbx^*$ versus the number of iterations. The exact solution $\bbx^*$ can be found in closed form fo the quadratic problem in \ref{eq_simulation_problem} and the relative error is then evaluated as
\begin{align}
\text{error}_{t} = \frac{1}{n} \sum_{i=1}^n \frac{\| \bbx_{i,t}  - \bbx^*\|^2}{\|\bbx^*\|^2}.
\label{eq_error_def}
\end{align}

Figure \ref{fig:4a} shows the convergence rates of all algorithms in the quadratic problem with small condition number $\eta=0$. Observe that PD-QN and D-ADMM converge substantially faster than all other methods, achieving an average error of $10^{-10}$ by iteration 100. The closest performing method in this simulation is the exact second order method ESOM, which doesn't quite reach an error of $10^{-8}$ after 400 iterations. While ESOM uses exact Hessian information in the primal update, it uses only a first order update in the dual domain, whereas PD-QN uses a quasi-Newton update in the dual domain. The decentralized ADMM method, on the other hand, is able to solve the augmented Lagraign of the quadratic objective easily in closed form, and thus experiences strong convergence rates with just first order information. For a larger condition number $\eta=1$, the corresponding convergence paths are given in Figure \ref{fig:4b}. Here, the performance of first order methods DA and D-ADMM degrades, as is commonly the case for first order methods in more ill-conditioned problems. The performance of the second order methods, as expected, degrade as well but still outperform the first order methods. PD-QN is the fastest converging method, reaching an error of $10^{-10}$ after 600 iterations. 

\begin{table}[h]
\begin{tabular}{ll}
       & Average Runtime (ms) \\\hline
PD-QN  & 2.2                  \\
DA     & 0.22                 \\
D-ADMM & 0.25                 \\
ESOM   & 0.87                 \\
D-BFGS & 6.7                 
\end{tabular}
\caption{Average runtimes per iteration for each method being compared.}
\label{tab_run}
\end{table}

The average runtime per iteration per method is reported in Table \ref{tab_run}. It can be seen that the quasi-Newton methods have higher runtime per runtime per iteration than the first order methods and ESOM. Both first order methods and ESOM will have very low complexity per iteration for quadratic problems such as the one considered here due to the fact that the augmented Lagrangian function can be solved in closed form and has a constant Hessian function. Generally speaking, higher order methods are more beneficial in scenarios in which communication complexity is of larger concern that computation complexity, as is often the case in settings such as sensor networks and distributed computing.

\begin{figure}
\centering
\includegraphics[width=.45\textwidth,height=.25\textheight,keepaspectratio]{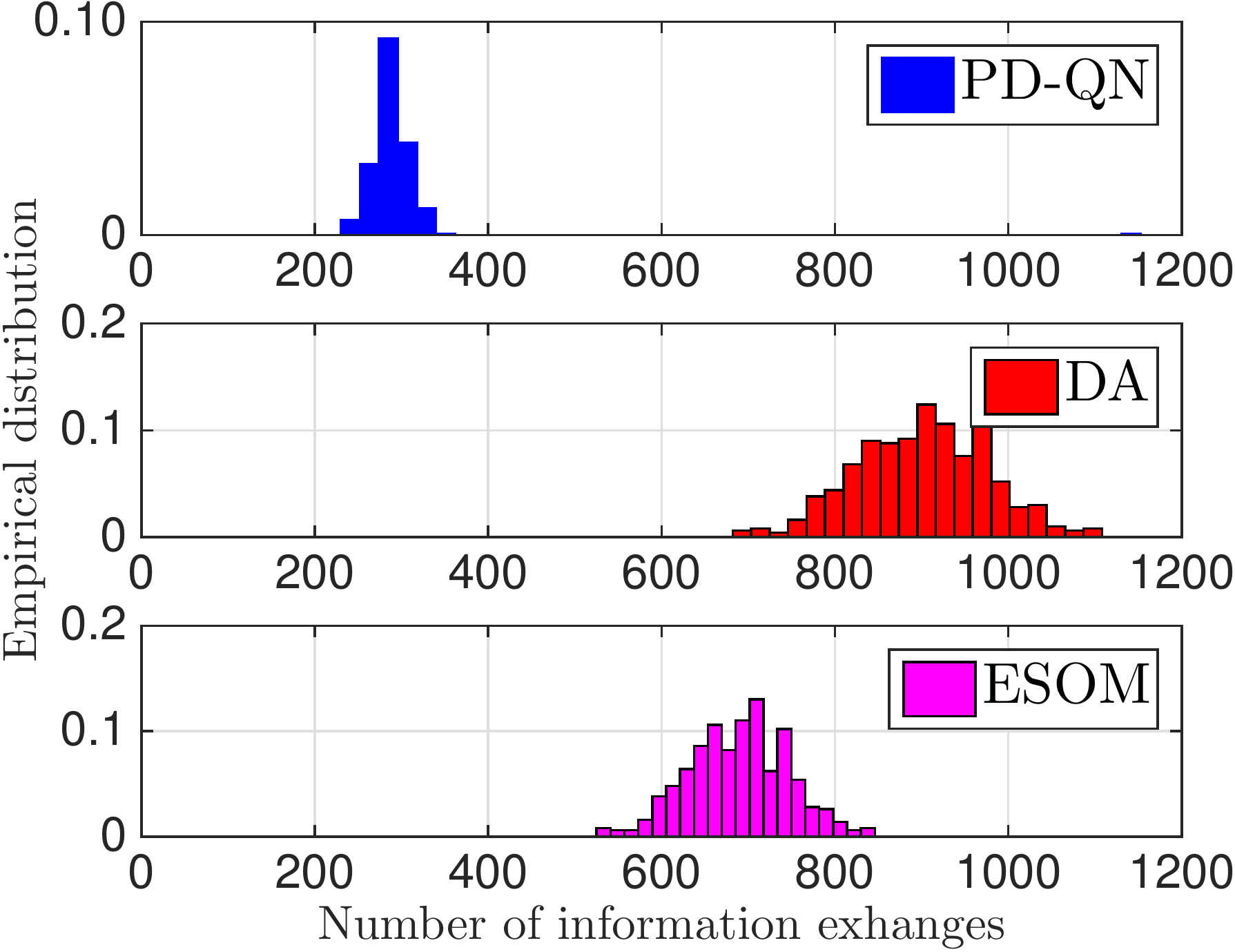}
\caption{Empirical distribution of number of information exchanges needed to reach error of $10^{-5}$ for PD-QN, DA, and ESOM for quadratic objective function with small condition number. The convergence gain of PD-QN is great enough that, even with the larger communication overhead, it outperforms the first and second order method when comparing information exchanges.}
\label{fig_primal_hist}
\end{figure}

\begin{figure}
\centering
\includegraphics[width=.45\textwidth,height=.2\textheight,keepaspectratio]{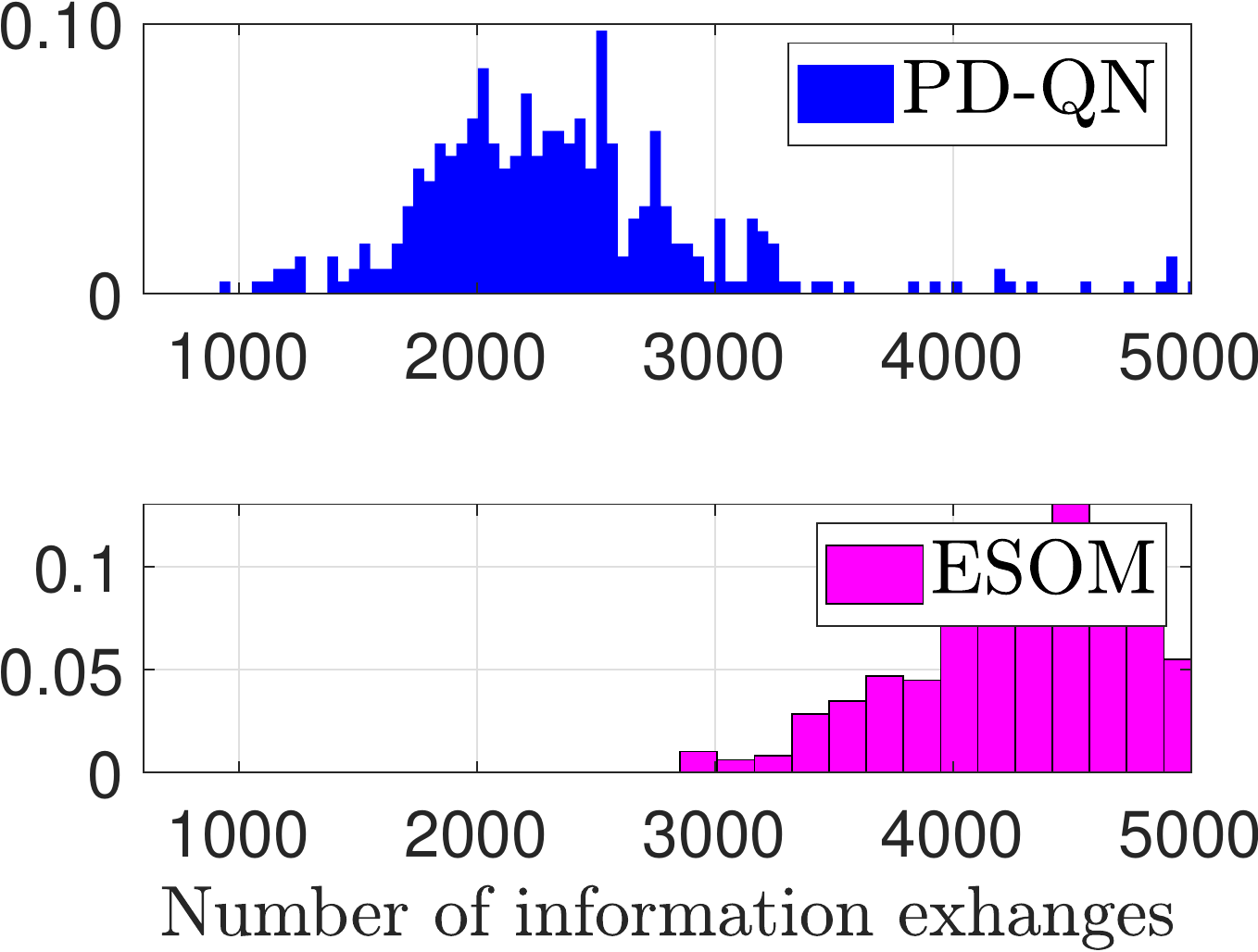}
\caption{Empirical distribution of number of information exchanges needed to reach error of $10^{-7}$ for PD-QN and ESOM for quadratic objective function with large condition number. PD-QN outperforms the second order ESOM method even with the additional communication overhead. }
\label{fig_primal_hist2}
\end{figure}

As previously discussed in Remark \ref{remark_comm}, the communication burden of PD-QN is larger than the alternatives, making the comparison in terms of number of iterations in Figure \ref{figure_simulation_quad} not entirely complete. To provide a more comprehensive comparison, we display in Figure \ref{fig_primal_hist} an empirical distribution of the performance of PD-QN, DA, and ESOM on the quadratic program with small condition number over 1000 different randomly generated experiments . In this case, however, we display the number of information exchanges needed to reach an error of $10^{-5}$. For this problem, observe that the PD-QN requires around 300 local information exchanges per node while DA and ESOM require around 700 and 900 exchanges, respectively. Indeed, the higher communication burden of PD-QN does not in this case outweigh the gain in number of iterations. For such a comparison for an ill-conditioned problem, we present in Figure \ref{fig_primal_hist2} the empirical distributions of communication exchanges required for the problem with large condition number for the quasi-Newton PD-QN method and the exact second order ESOM method. Overall, even with the additional communication overhead the PD-QN outperforms the existing primal-dual second order alternative method by roughly half the amount of communications. This highlights the additional benefit of the quasi-Newton update in the dual domain in the PD-QN method in problems with larger condition numbers.


\subsection{Logistic regression}
We perform additional simulations to evaluate the performance of PD-QN on a more complex objective function with varying condition number, namely distributed logistic regression problem that is very common in machine learning. We seek to learn a linear classifier $\bbx$ to predict the binary label of a data point $v_j \in \{-1,1\}$ given a feature vector $\bbu_j \in \reals^p$. For a set of training samples, we compute the likelihood of a label given a feature vector as $P(v=1 | \bbu) = 1/(1+\exp(-\bbu^T\bbx))$ and find the classifier $\bbx$ that maximizes the log likelihood over all samples. In the distributed setting, we assume that the training set is large and distributed amongst $n$ nodes, each holding r$q$ samples. Each node $i$ then has access to an objective function that is the loss over its training samples $\{\bbu_{il}\}_{l=1}^{q_i}$ and $\{v_{il}\}_{l=1}^{q_i}$. The aggregate objective function can be defined as 
\begin{align}
f(\bbx) := \frac{\lambda}{2}\| \bbx\|^2 + \sum_{i=1}^{n} \sum_{l=1}^{q_i} \log[ 1 + \exp(-v_{il}\bbu_{il}^T\bbx)],
\label{eq_logistic_problem}
\end{align}
where the first term is a regularization term used to reduce overfitting and is parametrized by $\lambda \geq 0$. 

\begin{figure*}[t]
\centering
	\begin{subfigure}[t]{.43\textwidth}
		\centering
		\includegraphics[height=\textheight,width=\textwidth,keepaspectratio]{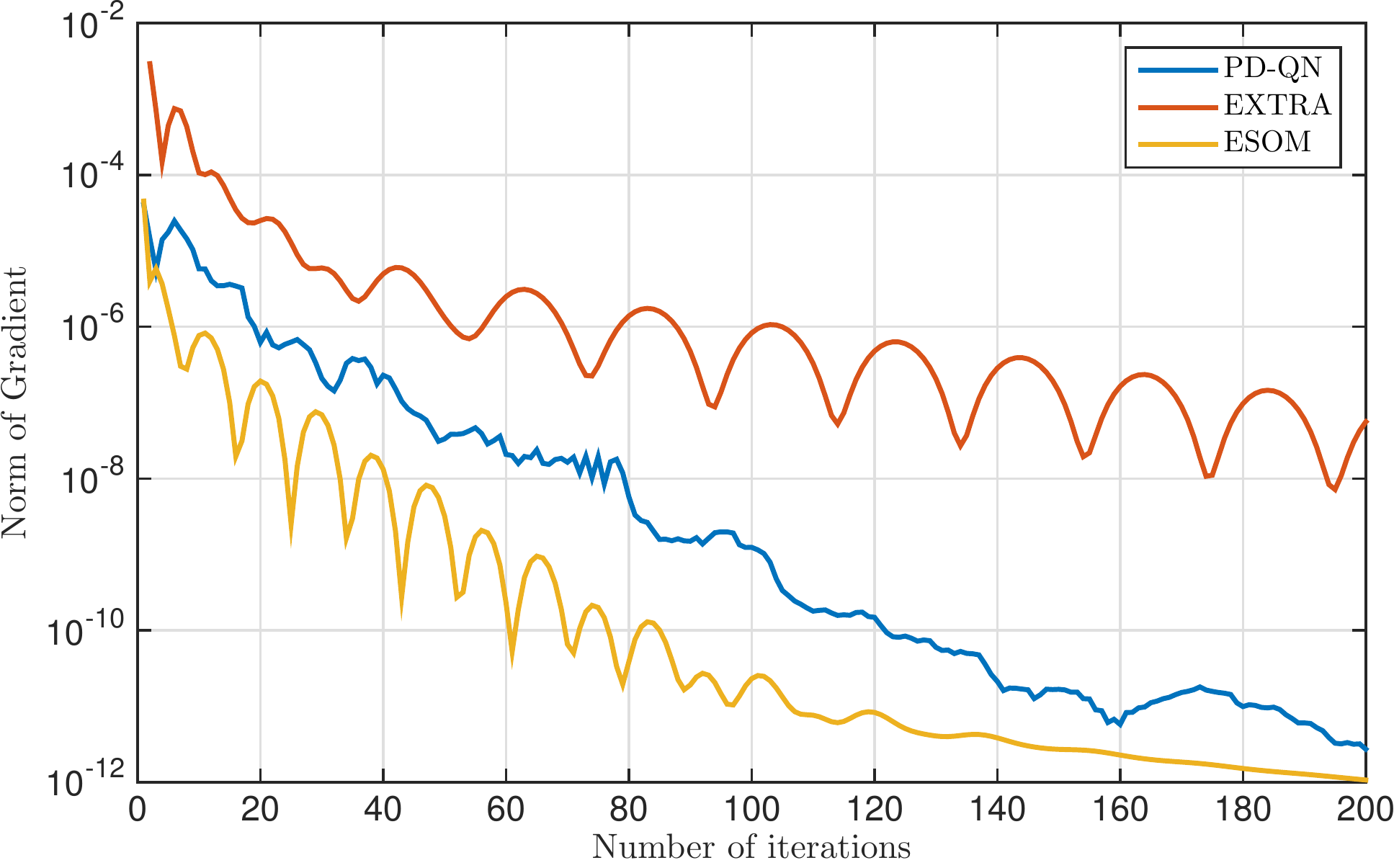}
		\caption{}\label{fig:5a}
	\end{subfigure}
	\begin{subfigure}[t]{.43\textwidth}
		\centering
		\includegraphics[height=\textheight,width=\textwidth,keepaspectratio]{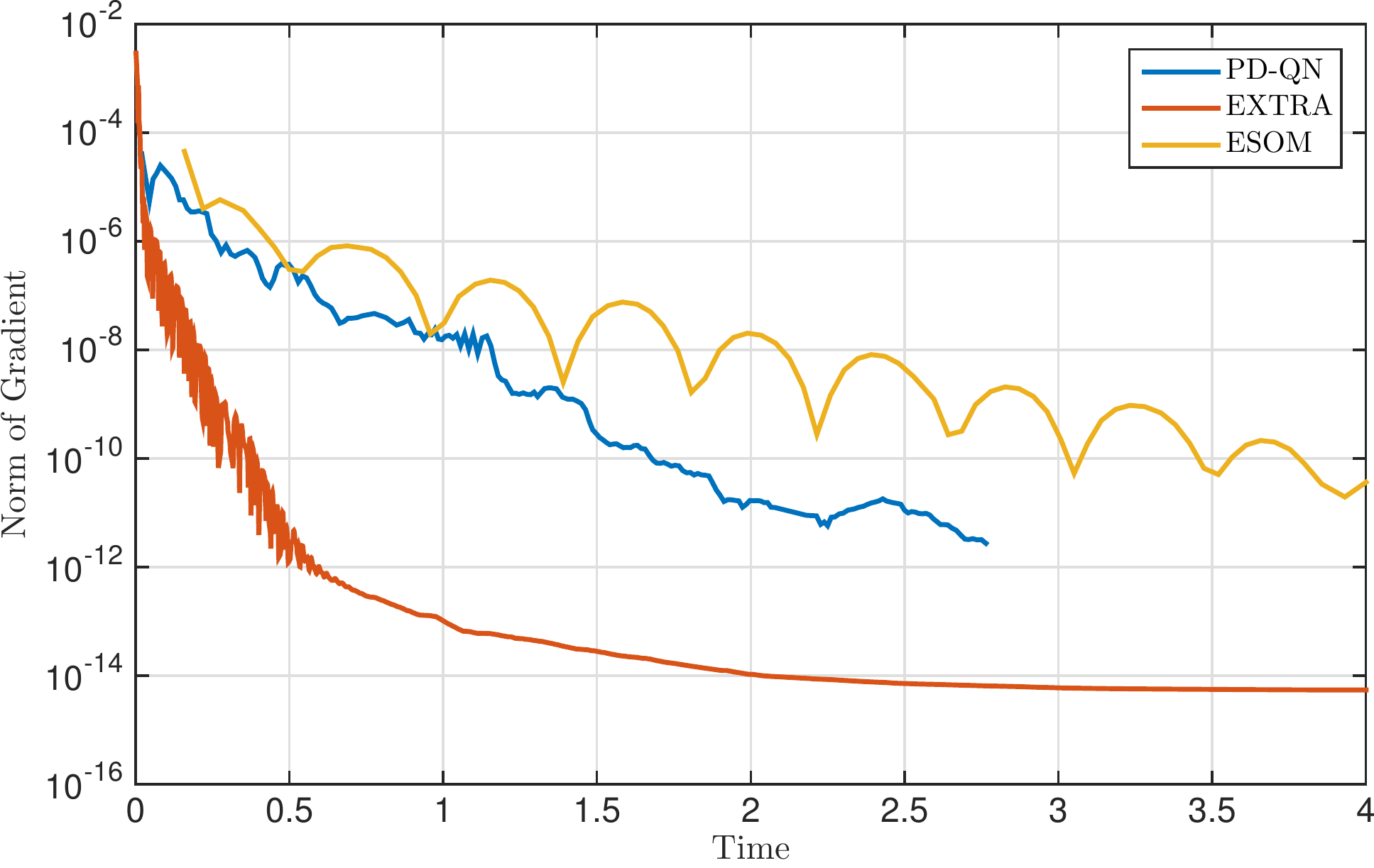}
		\caption{}\label{fig:5b}
	\end{subfigure} 	
	\caption{Convergence paths for exact distributed methods on the logistic regression problem, compared by (a) number of iterations and (b) runtime. The PD-QN method can be observed to split the difference between first order EXTRA method and second order ESOM method in terms of both metrics. }\label{figure_simulation_log}
\end{figure*}

For our simulations we generate an artificial dataset of feature vectors $\bbu_{il}$ with label $v_{il}=1$ from a normal distribution with mean $\mu$ and standard deviation $\sigma_{+}$, and with label $v_{il}=-1$ from a normal distribution with mean $-\mu$ and standard deviation $\sigma_{-}$. Each node $i$ receives $q_i=100$ samples and the regularization parameters is fixed to be $\lambda= 10^{-4}$. The feature vector parameters are set as $\mu=3$ and $\sigma_{+}=\sigma_{-}=1$ to make the data linearly separable. The other parameters we set the same as in earlier simulations, i.e. $n=20$ nodes connected in $d=4$-regular cycle with $p=4$. The PD-QN regularization parameters are chosen as $\Gamma = \gamma = 10^{-1}$.

For the logistic regression simulations, the form of the objective function in \eqref{eq_logistic_problem} does not permit an easily computable primal minimizer as in the case of the quadratic problem in \eqref{eq_simulation_problem}. Therefore, performing dual methods such as DA, ADMM, and D-BFGS require internal optimization problems at each iteration and are generally infeasible in this type of setting. Therefore, for these simulations we compare the performance of PD-QN only against the ``primal-dual'' type methods: first order EXTRA and second-order ESOM.  

The resulting convergence paths are shown in Figure \ref{figure_simulation_log}. The figure on the left compares the convergence in terms of number of iterations while the figure on the right compares convergence in terms of runtime. The results demonstrate a case in which the approximate second order PD-QN method splits the difference between the first and second order methods in terms of both metrics. When compared in terms of number of iterations, PD-QN performs similarly but slightly worse to ESOM, which uses exact Hessian information, while both methods outperform the first order EXTRA method. This is result reflects the fact the logistic regression problem has a more complex Hessian that is not as easily approximated with quasi-Newton methods. The ESOM method therefore benefits for using exact second order information, but at higher computational cost.  Indeed, the runtime comparison shows that EXTRA outperforms both methods because it has very low computational cost. Additionally, PD-QN outperforms ESOM in this manner due to the fact that computing exact Hessians for the logistic regression problem can be very costly, which thus makes it preferable to approximate the curvature information using the quasi-Newton updates of PD-QN.  PD-QN can thus be seen here to balance the tradeoff between the iteration complexity benefits of second order computation seen in ESOM with the runtime benefits of first order methods seen in EXTRA.

\subsection{Discussion}

In this section, we compare the performance of PD-QN against a number of popular decentralized consensus methods that use either first or second order information for both a distributed linear least squares, or quadratic, problem and a distributed logistic regression problem. In these results we can observe a number of tradeoffs between iteration complexity, computational complexity, and communication complexity. The PD-QN method outperforms the other methods because it estimates a second-order step in both the primal and dual variable updates, whereas the other methods perform first order updates in either one or both of the primal and dual variables updates. The use of second order information allows for speed up in both the maximization and minimization of the augmented Lagrangian function. The benefit of the second order update can be observed to be even greater when problems with larger condition numbers are considered. While the use of second order updates in PD-QN leads to higher computational complexity, we observe that the iteration complexity benefits generally outweigh the computational complexity issues in comparison to other first and second order decentralized methods.

%
\section{Conclusion} \label{sec_conclusion}
We considered the problem of decentralized consensus optimization, in which nodes sought to minimize an aggregate cost function while only being aware of a local strictly convex component. The problem was solved in the dual domain through the introduction of PD-QN as a decentralized quasi-Newton method. In PD-QN, each node approximates the curvature of its local cost function and its neighboring nodes to correct its descent direction. Analytical and numerical results were established showing its convergence and improvement over existing consensus methods, respectively.

\begin{appendices}

%
\section{Proof of Lemma \ref{lemma_pdqn_error}}\label{app_lemma_pdqn_error}

The details here follow closely those of a similar lemma in \cite{mokhtari2016decentralized}. Consider the primal and dual updates of PD-QN in \eqref{eq_pdqn_final1} and \eqref{eq_pdqn_final2}. 
To prove the result in \eqref{opt_res_PD-QN2}, we being by recalling the primal gradient $\bbg_{t} = \nabla f(\bbx_{t}) + \bby_{t} + \alpha (\bbI-\bbZ)\bbx_{t}$ and rearrange terms in \eqref{eq_pdqn_final1} to obtain
\begin{align}\label{pdqn_gen_proof_001}
\nabla f(\bbx_{t})+& \bby_{t}+{\alpha}(\bbI-\bbZ)\bbx_{t} +\bbG_{t,K}(\bbx_{t+1}-\bbx_{t})=\bb0.
\end{align}
Define $\bbL_{t}$ to be the Hessian of the augmented Lagrangian $\bbL_{t}:=\nabla^2_{\bbx \bbx} \ccalL_{\alpha}(\bbx_{t},\bby_{t}) = \nabla^2 f(\bbx_{t}) + \alpha (\bbI-\bbZ)$. We add and subtract the term $\bbL_{t}(\bbx_{t+1}-\bbx_{t})$ to \eqref{pdqn_gen_proof_001} to obtain
\begin{align}\label{pdqn_gen_proof_002}
&\nabla f(\bbx_{t})+\nabla^2f(\bbx_{t})(\bbx_{t+1}-\bbx_{t})+  \bby_{t}
\\&
+{\alpha}(\bbI-\bbZ)\bbx_{t+1}+(\bbG_{t,K}- \bbL_{t})(\bbx_{t+1}-\bbx_{t})=\bb0.\nonumber
\end{align}
Now using the definition of the error vector $\bbsigma_{t}$ in \eqref{esom_error_vec} we can rewrite \eqref{pdqn_gen_proof_002} as
\begin{align}\label{pdqn_gen_proof_003}
&\nabla f(\bbx_{t+1})+\bby_{t}
+{\alpha}(\bbI-\bbZ)\bbx_{t+1}
+\bbsigma_{t} \\
&\qquad +\alpha\bbH_{t}^{-1}(\bbI-\bbZ)(\bbx_{t+1}-\bbx^*) =\bb0, \nonumber
\end{align}
where we use the fact that $\bbL_{t} - \nabla^2f(\bbx_{t}) = \alpha(\bbI-\bbZ)$ in the definition of $\bbsigma_{t}$ in \eqref{esom_error_vec}.
To prove the claim in \eqref{opt_res_PD-QN2} from \eqref{pdqn_gen_proof_003}, first consider that one of the KKT conditions of the optimization problem in \eqref{eq_primal_problem} is  
\begin{equation}\label{KKT_condition}
\nabla f(\bbx^*) + (\bbI-\bbZ)^{1/2}\bbnu^*=\nabla f(\bbx^*) +\bby^*=\bb0.
\end{equation}
Subtracting the equality in \eqref{KKT_condition} from \eqref{pdqn_gen_proof_003} yields
\begin{align}\label{important_result_100}
&\nabla f(\bbx_{t+1})-\nabla f(\bbx^*) +\bby_{t}-\bby^* \nonumber\\
&\quad+\alpha\bbH_{t}^{-1}(\bbI-\bbZ) (\bbx_{t+1}-\bbx^*) +\bbsigma_{t} =\bb0.
\end{align}
Furthermore, by using the dual variable update in \eqref{eq_pdqn_final2} along with the consensus constraint $\alpha(\bbI-\bbZ)\bbx^*=\bb0$, we can additionally claim that
\begin{equation}\label{important_result_200}
\bby_{t} =\bby_{t+1}-\alpha \bbH^{-1}_{t} (\bbI-\bbZ)(\bbx_{t+1}-\bbx^*).
\end{equation}
Substituting $\bby_{t}$ in \eqref{important_result_100} by the expression in the right hand side of \eqref{important_result_200} leads to the claim in \eqref{opt_res_PD-QN2}.

%
\section{Proof of Theorem \ref{thm:esom_linear_convg}}\label{app_thm_esom_linear_convg}

The details here are again adapted from a similar result in \cite{mokhtari2016decentralized}, but modified to consider the quasi-Newton primal and dual updates present in PD-QN. To prove this result, we begin by applying a well-known lower bound for the inner product $(\bbx_{t+1}-\bbx_{t})^T(\nabla f(\bbx_{t+1})-\nabla f(\bbx^*))$ that incorporates both strong convexity constant $\mu$ and the Lipschitz constant of the gradients $L$. This inequality can be written as
\begin{align}\label{proof_lin_x}
&\frac{1}{\mu + L}(\mu L\|\bbx_{t+1}-\bbx^*\|^2 + \|\nabla f(\bbx_{t+1})-\nabla f(\bbx^*)\|^2) \nonumber \\
&\qquad\leq (\bbx_{t+1}-\bbx^*)^T(\nabla f(\bbx_{t+1})-\nabla f(\bbx^*)).
\end{align}
 The result in \eqref{opt_res_PD-QN2} gives us an expression for the difference $\nabla f(\bbx_{t+1})-\nabla f(\bbx^*)$ by rearranging terms as
\begin{align}\label{proof_lin_pdqn_001}
\nabla f(\bbx_{t+1})-\nabla f(\bbx^*) &=  -\bby_{t+1}+\bby^* - \bbsigma_{t}.
\end{align}
We now substitute \eqref{proof_lin_pdqn_001} into  \eqref{proof_lin_x} and multiply both sides of the inequality by $2\alpha$ to obtain
\begin{align}\label{proof_lin_pdqn_002}
&
\frac{2\alpha \mu L}{\mu+L}\|\bbx_{t+1}-\bbx^*\|^2+\frac{2\alpha}{\mu+L}\|\df(\bbx_{t+1})-\df(\bbx^*)\|^2\nonumber\\
& \leq -2\alpha(\bbx_{t+1}-\bbx^*)^T(\bby_{t+1}-\bby^*) -2\alpha(\bbx_{t+1}-\bbx^*)^T\bbsigma_{t}.
\end{align}
To proceed from here, we will substitute the modified dual variable $\bby_{t}$ back to the original dual variable $\bbnu_{t}$ by recalling the transformation $\bby_{t} = (\bbI-\bbZ)^{1/2}\bbnu_{t}$. First, we rewrite \eqref{proof_lin_pdqn_002} as
\begin{align}\label{proof_lin_pdqn_003}
&
\frac{2\alpha \mu L}{\mu+L}\|\bbx_{t+1}-\bbx^*\|^2+\frac{2\alpha}{\mu+L}\|\df(\bbx_{t+1})-\df(\bbx^*)\|^2\nonumber\\
& \leq -2\alpha(\bbx_{t+1}-\bbx^*)^T(\bbI-\bbZ)^{-1/2}(\bbnu_{t+1}-\bbnu^*) \nonumber \\
& \qquad -2\alpha(\bbx_{t+1}-\bbx^*)^T\bbsigma_{t}.
\end{align}
Now we can derive a similar expression as in \eqref{important_result_200} with respect to $\bbnu_{t}$. Again using the fact that $(\bbI-\bbZ)^{-1/2}\bbx^* = \bb0$, we can add this term to the dual update in \eqref{eq_pdqn_final2} in terms of $\bbnu_{t}$, rearrange terms to obtain that
\begin{equation}\label{proof_lin_pdqn_004}
\alpha(\bbI-\bbZ)^{-1/2}(\bbx_{t+1}-\bbx^*) = \bbH_{t}(\bbnu_{t+1}-\bbnu_{t}).
\end{equation}
Now we can substitute \eqref{proof_lin_pdqn_004} back into \eqref{proof_lin_pdqn_003} to obtain
\begin{align}\label{proof_lin_pdqn_005}
&
\frac{2\alpha \mu L}{\mu+L}\|\bbx_{t+1}-\bbx^*\|^2+\frac{2\alpha}{\mu+L}\|\df(\bbx_{t+1})-\df(\bbx^*)\|^2\nonumber\\
& \leq -2(\bbnu_{t+1}-\bbnu_{t})^T\bbH_{t} (\bbnu_{t+1}-\bbnu^*)  -2\alpha(\bbx_{t+1}-\bbx^*)^T\bbsigma_{t}.
\end{align}
From here, we wish to write the primal and dual terms in \eqref{proof_005} in terms of a combined variable $\bbz_{t} := [\bbx_{t}; \bbnu_{t}]$. To that end, we first decompose the error term as $\bbsigma_{t}$ as $\bbsigma_{t} = \hbsigma_{t} + \left[\bbG_{t,K}-\alpha(\bbI-\bbZ)\right](\bbx_{t+1}-\bbx_{t})$, where $\hbsigma_{t} := \nabla f(\bbx_{t})-\nabla f(\bbx_{t+1}) + \alpha\bbH_{t}^{-1}(\bbI-\bbZ)(\bbx_{t+1}-\bbx^*)$. Now we rewrite \eqref{proof_lin_pdqn_005} as
\begin{align}\label{proof_0051}
&
\frac{2\alpha \mu L}{\mu+L}\|\bbx_{t+1}-\bbx^*\|^2+\frac{2\alpha}{\mu+L}\|\df(\bbx_{t+1})-\df(\bbx^*)\|^2\nonumber\\
& \leq 
 -2(\bbnu_{t+1}-\bbnu_{t})^T\bbH_{t} (\bbnu_{t+1}-\bbnu^*) \nonumber\\  &  -2\alpha(\bbx_{t+1}-\bbx^*)^T \bbR_{t}(\bbx_{t+1}-\bbx_{t})  -2\alpha(\bbx_{t+1}-\bbx^*)^T\hbsigma_{t}, 
\end{align}
where we define $\bbR_{t} := \left[\bbG_{t,K}-\alpha(\bbI-\bbZ)\right]$ for notational simplicity.
We transform the vector difference products on the right hand side of \eqref{proof_0051} using the distributive property. For any vectors $\bba$, $\bbb$, and $\bbc$ we can write $
2(\bba-\bbb)^T(\bba-\bbc)=\|\bba-\bbb\|^2+\|\bba-\bbc\|^2-\|\bbb-\bbc\|^2$. By applying this substitution into \eqref{proof_0051} for the first two terms on the right hand side, we have that
\begin{align}\label{proof_005}
&
\frac{2\alpha \mu L}{\mu+L}\|\bbx_{t+1}-\bbx^*\|^2+\frac{2\alpha}{\mu+L}\|\df(\bbx_{t+1})-\df(\bbx^*)\|^2\nonumber\\
& \leq 
\left( \|\bbnu_{t}-\bbnu^*\|^2_{\bbH_t} -\|\bbnu_{t+1}-\bbnu_{t}\|^2_{\bbH_t}-  \|\bbnu_{t+1}  -\bbnu^*\|^2_{\bbH_t} \right)  \nonumber \\ 
&+ \alpha\left( \|\bbx_{t}-\bbx^*\|^2_{\bbR_t} -\|\bbx_{t+1}-\bbx_{t}\|^2_{\bbR_t}-  \|\bbx_{t+1}  -\bbx^*\|^2_{\bbR_t} \right)  \nonumber \\
&\qquad -2\alpha(\bbx_{t+1}-\bbx^*)^T\hbsigma_{t}.
\end{align}

We make an additional substitution to upper bound \eqref{proof_005}. Observe from \eqref{proof_lin_pdqn_004} that we have the equivalence $\|\bbnu_{t+1}-\bbnu_{t}\|_{\bbH_t^2}^2 = \|\bbx_{t+1}-\bbx^*\|_{\alpha^2(\bbI-\bbZ)}^2$. We may further use the upper and lower bounds of the eigenvalues of $\bbH_t$ in \eqref{eq_prop_eigen_bounds_dual} to show that $\| \bby \|^2_{\bbH_t^2} \leq \| \bby \|^2_{\bbH_t}$ for any vector $\bby$. In particular, we may consider the upper bound $\| \bby \|^2_{\bbH_t^2} \leq \| \bby \|^2 / \Gamma^2$  and the lower bound $\| \bby \|^2/P \leq \| \bby \|^2_{\bbH_t}$ and show that the former upper bound is lower than the latter lower bound, i.e.
\begin{align}
\frac{1}{\Gamma^2} \| \bby\|^2 &\leq \frac{1}{P} \|\bby\|^2 =  \left(\frac{1}{\Gamma + n/\gamma}\right) \|\bby\|^2 \nonumber \\
\frac{1}{\Gamma} - 1 &\leq \frac{\Gamma \gamma}{n}. \label{eq_proof_eq}
\end{align}
The inequality in \eqref{eq_proof_eq} holds when $\Gamma \leq 1$, as the left hand side will be negative while the right hand side is positive. 

We proceed with the main proof by introduce the matrix $\bbJ_t := \diag(\alpha \bbR_t,  \bbH_t)$ and the combined vector $\bbz_t := [\bbx_t; \bbnu_t]$ to combine the terms in \eqref{proof_005}. Furthermore, we can substitute from \eqref{proof_lin_pdqn_004} and the logic in \eqref{eq_proof_eq} that $-\|\bbnu_{t+1}-\bbnu_{t}\|_{\bbH_t} \leq  -\|\bbx_{t+1}-\bbx^*\|_{\alpha^2(\bbI-\bbZ)}^2$ and rearrange terms in \eqref{proof_005} to obtain
\begin{align}\label{proof_007}
&
\|\bbz_t-\bbz^*\|_{\bbJ_t}^2 -\|\bbz_{t+1}-\bbz^*\|_{\bbJ_t}^2
\\
&\quad\geq 
\frac{2\alpha}{\mu+L}\|\df(\bbx_{t+1})-\df(\bbx^*)\|^2
+\alpha \Sigma\|\bbx_{t+1}-\bbx_{t}\|^2
\nonumber\\
&\qquad +\|\bbx_{t+1}-\bbx^*\|^2_{\frac{2\alpha \mu L}{\mu+L}\bbI+\alpha^2(\bbI-\bbZ)} +2\alpha(\bbx_{t+1}-\bbx^*)^T\hbsigma_{t}.\nonumber
\end{align}
Note that the inner product $2(\bbx_{t+1}-\bbx^*)^T\hbsigma_{t}$ is bounded below by $-(1/\zeta)\|\bbx_{t+1}-\bbx^*\|^2-\zeta \|\hbsigma_{t}\|^2$ for any positive constant $\zeta>0$. Thus, the lower bound in \eqref{proof_007} can be updated as
\begin{align}\label{proof_lin_pdqn_004b}
&\|\bbz_{t}-\bbz^*\|_{\bbJ_t}^2-\|\bbz_{t+1}-\bbz^*\|_{\bbJ_t}^2\nonumber\\
& \geq
\|\bbx_{t+1}-\bbx^*\|_{(\frac{2\alpha \mu L}{\mu+L}-\frac{\alpha}{\zeta})\bbI+\alpha^2(\bbI-\bbZ)}^2+\alpha \Sigma\|\bbx_{t+1}-\bbx_{t}\|^2
 \nonumber\\
&\ +\frac{2\alpha}{\mu+L}\|\df(\bbx_{t+1})-\df(\bbx^*)\|^2-\alpha \zeta\|\hbsigma_{t}\|^2.
\end{align}
The linear convergence result in \eqref{pdqn_lin_convg} is equivalent to establishing the inequality $\|\bbz_{t}-\bbz^*\|_{\bbJ_t}^2-\|\bbz_{t+1}-\bbz^*\|_{\bbJ_t}^2 \geq \kappa_{t}\|\bbz_{t+1}-\bbz^*\|_{\bbJ}^2$. In other words, we may lower bound the right hand side of \eqref{proof_lin_pdqn_004b}  by $\kappa_{t}\|\bbz_{t+1}-\bbz^*\|_{\bbJ_t}^2$, i.e., 
\begin{align}\label{proof_lin_pdqn_005b}
&\kappa_{t}\|\bbnu_{t+1}-\bbnu^*\|_{\bbH_t}^2+\kappa_{t} \alpha \|\bbx_{t+1}-\bbx^*\|_{\bbR_t}^2\nonumber\\
& \leq
\|\bbx_{t+1}-\bbx^*\|_{(\frac{2\alpha \mu L}{\mu+L}-\frac{\alpha}{\zeta})\bbI+\alpha^2(\bbI-\bbZ)}^2+\alpha \|\bbx_{t+1}-\bbx_{t}\|_{\bbR_t}^2
 \nonumber\\
&\ +\frac{2\alpha}{\mu+L}\|\df(\bbx_{t+1})-\df(\bbx^*)\|^2-\alpha \zeta\|\hbsigma_{t}\|^2.
\end{align}
We will determine values of $\kappa_t$ for which the inequality in \eqref{proof_lin_pdqn_005b} is satisfied. A sufficient condition can be formulated by lower bounding the left hand side of the inequality in \eqref{proof_lin_pdqn_005b}. We may lower bound the terms $\|\bbnu_{t+1}-\bbnu^*\|^2_{\bbH_t} + \alpha \|\bbx_{t+1}-\bbx^*\|^2_{\bbR_t}$ by using the lower eigenvalue bound of $\bbH_t$---namely $P^{-1}$ from \eqref{eq_prop_eigen_bounds_dual} and the lower eigenvalue bound of $\bbR_t$---namely $\Lambda^{-1} - 2\alpha(1-\delta)$. Observe that this lower bound is obtained by combining the lower eigenvalue bound of $\bbG_{t,K}$ in \eqref{eq_prop_eigen_bounds_primal} with that of $-\alpha(\bbI-\bbZ)$ that can be found in, e.g. \cite[Proposition 1]{mokhtari2017network}.
From these bounds, we obtain
\begin{align}\label{proof_lin_pdqn_005c}
&\frac{\kappa_{t}}{P}\|\bbnu_{t+1}-\bbnu^*\|^2+\kappa_{t} \alpha \Sigma \|\bbx_{t+1}-\bbx^*\|^2\nonumber\\
& \leq
\|\bbx_{t+1}-\bbx^*\|_{(\frac{2\alpha \mu L}{\mu+L}-\frac{\alpha}{\zeta})\bbI+\alpha^2(\bbI-\bbZ)}^2+\alpha \|\bbx_{t+1}-\bbx_{t}\|^2_{\bbR_t}
 \nonumber\\
&\ +\frac{2\alpha}{\mu+L}\|\df(\bbx_{t+1})-\df(\bbx^*)\|^2-\alpha \zeta\|\hbsigma_{t}\|^2,
\end{align}
where we define $\Sigma := \Lambda^{-1} - 2\alpha(1-\delta)$ for notational convenience.

From here, we first find an upper bound for the term $\|\bbnu_{t+1}-\bbnu^*\|^2$. Consider the relation for $\bbnu_{t+1}-\bbnu^*$ derived from substituting $\bby_t = (\bbI-\bbZ)^{1/2}\bbnu_t$ into \eqref{opt_res_PD-QN2}. Further consider for any $\beta, \phi > 1$ and arbitrary vectors $\bba$, $\bbb$, $\bbc$, we can write that $(1-\beta^{-1})(1-\phi^{-1})\|\bbc^2\| \leq \|\bba + \bbb+\bbc\|^2+(\beta-1)\|\bba\|^2+(\phi-1)(1-\beta^{-1})\|\bbb\|^2$. Combining these we have that
\begin{align}\label{proof_lin_pdqn_006}
\|\bbnu_{t+1}-\bbnu^*\|^2
&\leq
\frac{\beta^2}{(\beta-1)\Gamma \hat{\delta}}\|\bbx_{t+1}-\bbx_{t}\|_{\bbR_{t}}^2+\frac{\beta\phi}{(\phi-1)\Gamma\hat{\delta}}\|\hbsigma_{t}\|^2
\nonumber\\
& \quad + \frac{\phi\beta}{\Gamma\hat{\delta}}\|\df(\bbx_{t+1})-\df(\bbx^*)\|^2 ,
\end{align}
where we define $\hat{\delta}$ as the smallest non-zero eigenvalue of $(\bbI-\bbZ)$ to simplify notation. By substituting the upper bound in \eqref{proof_lin_pdqn_006} for the squared norm $\|\bbnu_{t+1}-\bbnu^*\|^2$ in \eqref{proof_lin_pdqn_005c} we obtain a sufficient condition for the result in \eqref{proof_lin_pdqn_005c} which is given by 
\begin{align}\label{proof_lin_pdqn_007}
&\kappa_{t} \alpha \Sigma \|\bbx_{t+1}-\bbx^*\|
+
\frac{\kappa_{t}\beta^2}{P (\beta-1)\Gamma\hat{\delta}}\|\bbx_{t+1}-\bbx_{t}\|_{\bbR_{t}}^2
\nonumber\\
&+ \frac{\kappa_{t}\phi\beta}{P \Gamma\hat{\delta}}\|\df(\bbx_{t+1})-\df(\bbx^*)\|^2
+\frac{\kappa_{t}\beta\phi\|\hbsigma_{t}\|^2 }{P (\phi-1)\Gamma\hat{\delta}}\nonumber\\
&\quad  \leq
\|\bbx_{t+1}-\bbx^*\|_{(\frac{2\alpha \mu L}{\mu+L}-\frac{\alpha}{\zeta})\bbI+\alpha^2(\bbI-\bbZ)}^2+\alpha \|\bbx_{t+1}-\bbx_{t}\|_{\bbR_t}^2
 \nonumber\\
&\qquad +\frac{2\alpha}{\mu+L}\|\df(\bbx_{t+1})-\df(\bbx^*)\|^2-\alpha \zeta\|\hbsigma_{t}\|^2.
\end{align}
It remains to bound norm of the error term $\|\hbsigma_{t}\|^2$. Consider the standard inequality for squared norms that comes from the application of the Cauchy-Schwartz inequality, i.e.
\begin{align}\label{eq_alo}
\| \hbsigma_{t}\|^2 &\leq 2\| \nabla f(\bbx_{t})-\nabla f(\bbx_{t+1}) \|^2 \\
&\qquad + 2\| \alpha \bbH^{-1}_{t}(\bbI - \bbZ) (\bbx_{t+1}-\bbx^*)\|^2. \nonumber
\end{align}
To bound the first term on the right hand side of \eqref{eq_alo}, we can use the Lipschitz continuity of the gradient implied from Assumption \ref{as_strongly_convex} to obtain
\begin{align}
2\| \nabla f(\bbx_{t})-\nabla f(\bbx_{t+1}) \| &\leq 2 L^2 \| \bbx_{t+1} -\bbx_{t}\|^2 \nonumber\\
& \leq \frac{2L^2}{\Sigma} \| \bbx_{t+1} -\bbx_{t}\|_{\bbR_t}^2.\label{eq_alo1}
\end{align}
The second term in \eqref{eq_alo} can subsequently be bounded using the Cauchy-Schwartz inequality along with the bound on the magnitude of $\bbH^{-1}_{t}$ from \eqref{eq_prop_eigen_bounds_dual} to obtain
\begin{align}\label{eq_alo3}
2 \| \alpha \bbH^{-1}_{t} (\bbI - \bbZ) (\bbx_{t+1}-\bbx^*)\|^2 \leq 4\alpha^2 P(1-\delta) \| \bbx_{t+1}-\bbx^*\|^2.
\end{align}
Combining the results of \eqref{eq_alo1}-\eqref{eq_alo3}, we obtain 
\begin{align}\label{alg4}
\| \hbsigma_{t}\|^2 &\leq 2L^2 \| \bbx_{t+1} -\bbx_{t}\|^2 + 4\alpha^2 P(1-\delta) \| \bbx_{t+1}-\bbx^*\|^2.
\end{align}
We proceed by substituting $\| \hbsigma_{t}\|^2$ in \eqref{proof_lin_pdqn_007} by the upper bound that is derived in \eqref{alg4}. After rearranging terms, we obtain
\begin{align}\label{proof_lin_pdqn_008}
&  
0 \leq
\|\bbx_{t+1}-\bbx^*\|_{(\frac{2\alpha \mu L}{\mu+L}-
\frac{\alpha}{\zeta}-\kappa_{t} \alpha \Sigma -\frac{2\kappa_{t}\beta\phi \alpha^2  }{P(\phi-1)\hat{\delta}} -4\alpha^3P(1-\delta)\zeta)\bbI+\alpha^2(\bbI-\bbZ)}^2\nonumber\\
&\ +\left(\frac{2\alpha}{\mu+L}-\frac{\kappa_{t}\phi\beta}{P\Gamma\hat{\delta}}\right)\|\df(\bbx_{t+1})-\df(\bbx^*)\|^2 \nonumber\\
&\
+\!\!\Bigg[\alpha \Sigma\!-\!\frac{\kappa_{t}\beta^2}{P(\beta-1)\hat{\delta}}
\!-\!\frac{2\kappa_{t}\beta\phi L^2  }{P(\phi-1)\hat{\delta}}\!-\!\frac{2\alpha \zeta L^2}{\Sigma}  \Bigg]\|\bbx_{t+1}-\bbx_{t}\|_{\bbR_t}^2.
\end{align}

The inequality that is established in \eqref{proof_lin_pdqn_008} provides a condition on $\kappa_t$ that, when satisfied, implies the inequality in \eqref{proof_lin_pdqn_007} holds. This is turn implies \eqref{proof_lin_pdqn_005b} and subsequently, and most importantly, the linear convergence statement in \eqref{pdqn_lin_convg}. We may satisfy the inequality in \eqref{proof_lin_pdqn_008} by ensuring that the coefficients of the three terms in \eqref{proof_lin_pdqn_008} are non-negative. This is to say that  \eqref{proof_lin_pdqn_008} holds if $\kappa_{t}$ satisfies 

\begin{align}\label{havij}
&\frac{2\alpha \mu L}{\mu+L}-\frac{\alpha}{\zeta}-\kappa_{t} \alpha \Sigma -\frac{2\kappa_{t}\beta\phi\alpha^2  }{P(\phi-1)\hat{\delta}} -2\alpha^3P\zeta \geq 0, \\
&\frac{2\alpha}{\mu+L} \geq \frac{\kappa_{t}\phi\beta}{P \Gamma \hat{\delta}} \nonumber \\
&\alpha \Sigma \geq\frac {\kappa_{t}\beta^2}{P(\beta-1)\hat{\delta}}+\frac{2\kappa_{t}\beta\phi\Gamma^2  }{P(\phi-1)\hat{\delta}}+\frac{2\alpha \zeta L^2}{\Sigma}.\nonumber
\end{align}
Observe that the expressions in \eqref{havij} are satisfied when $\kappa_t$ is chosen as in the statement of Theorem \ref{thm:esom_linear_convg}, in which case the claim in \eqref{pdqn_lin_convg} holds.

\end{appendices}

\urlstyle{same}
\bibliographystyle{IEEEtran}
\bibliography{bmc_article}

\end{document}